\documentclass[11pt,a4paper]{article}
\ifdefined\pdfsuppresswarningpagegroup
\fi
\usepackage[T1]{fontenc}
\usepackage[utf8]{inputenc}
\usepackage[english]{babel}
\usepackage{amsmath,amssymb,amsthm,mathtools}
\usepackage{graphicx}
\usepackage{xcolor}
\usepackage[hidelinks]{hyperref}
\usepackage{caption}
\usepackage{booktabs}
\usepackage{tikz}
\mathtoolsset{showonlyrefs}

\newtheorem{theorem}{Theorem}[section]
\newtheorem{lemma}[theorem]{Lemma}
\newtheorem{prp}[theorem]{Proposition}
\newtheorem{assumption}[theorem]{Assumption}
\newtheorem{remark}[theorem]{Remark}

\newcommand{\R}{\mathbb{R}}
\newcommand{\N}{\mathbb{N}}

\title{A rescaling algorithm for finite-time quenching in a semilinear parabolic equation}
\author{Maissâ Boughrara \\
Fakultät für Mathematik, Universität Bielefeld,\\ D-33501 Bielefeld, Germany\vspace{1em}\\
Takiko Sasaki\\
Department of Mathematical Engineering,\\ Faculty of Engineering,
Musashino University,\\
Tokyo 135-8181, Japan
}
\date{}

\definecolor{revisionred}{HTML}{D20000}
\definecolor{revisionblue}{HTML}{0046BE}
\newcommand{\revred}[1]{{\color{revisionred}#1}}
\newcommand{\revblue}[1]{{\color{revisionblue}#1}}
\newtagform{blacktags}[\textcolor{black}]{\textcolor{black}{(}}{\textcolor{black}{)}}
\usetagform{blacktags}
\DeclareRobustCommand{\eqref}[1]{\noeqref{#1}\textup{(\ref{#1})}}
\let\restoredqedsymbol\qedsymbol
\renewcommand{\qedsymbol}{{\color{black}\restoredqedsymbol}}

\definecolor{maissared}{HTML}{D20000}
\definecolor{revisionred}{HTML}{000000}
\definecolor{revisionblue}{HTML}{000000}

\begin{document}
\maketitle

\begin{abstract}
In this work, we study numerical approximations of parabolic equations whose solutions exhibit finite-time quenching.  After introducing the change of variable $u=1-U$, we focus on the normalized one-dimensional singular problem
$$
    u_t=u_{xx}-u^{-r},\qquad u(\pm1,t)=1,
$$
and adapt the Berger and Kohn rescaling idea to the quenching setting.  We first analyze the associated ODE model and derive first-order estimates for the numerical quenching time.  For the PDE, the algorithm uses a nodal stopping time, linear interpolation in time, and local cubic reconstruction in space at every rescaling step. The convergence analysis is stated for any fixed finite number of rescaling stages and explicitly includes the nonzero initial error inherited from the previous numerical stage, the time interpolation error, and the cubic spatial reconstruction error.  The convergence of the full numerical quenching time is conditional on a uniform tail comparison assumption, which is also checked a posteriori in the numerical experiments.  The computations are performed in an aligned-grid setting consistent with the hypotheses of the analysis.
\end{abstract}

\section{Introduction}

We consider the following parabolic equation:
\begin{equation}\label{EQ}
\begin{cases}
    &\partial_t U = \Delta U    +
\displaystyle \frac{ a}{  (1 - U)^r }
,\ x\in\Omega, t>0,\\
&U(x,t)=0, t>0, x\in \partial \Omega,\\
&U(x,0)=U_0(x), x\in \Omega,
\end{cases}
\end{equation}
where $a,r>0$ and $\Omega\subset \R^N$. We are interested in quenching solutions for equation \eqref{EQ}. Quenching refers to the phenomenon that the distance to the singular state vanishes in finite time. In the present variables this means
\begin{equation}\label{eq:touchdown-U}
   \lim_{t\rightarrow T}\min_{x\in \overline\Omega}(1-U(x,t))=0
\end{equation}
for some $T<+\infty$. This type of equation has been studied extensively; see, for example, \cite{LM1980, AW1976, DZ2019}. In particular, for the existence of quenching solutions for equation \eqref{EQ}, we refer to the work of Levine \cite{L1989} in the one-dimensional case with $ r \geq 1 $ and $N=1$. In this work, Levine constructed the set of stationary solutions and analyzed their stability. One of the key results is as follows: there exists a critical value $ \varepsilon(r) > 0 $ such that for $ a > \varepsilon(r) $, all solutions exhibit quenching behavior. \revred{Moreover, we say that $x_0\in \Omega$ is a touch-down point if and only if there exists a sequence $(x_n,t_n)\in \Omega\times (0,T)$ such that}

$$(x_n, t_n) \rightarrow (x_0,T) \text{ and } U(x_n, t_n)\rightarrow 1, \text{ as } n\rightarrow +\infty.$$
For the study of quenching points, we cite \cite{GH2018, GS2015, N2011}. Acker and Kawohl \cite{AK1989} studied singular quenching points and proved the single-point quenching property for a class of radial initial conditions in higher dimensions. \revred{Moreover, Merle and Zaag \cite{MZ1997Nonlinearity} constructed a solution of \eqref{EQ}, in the normalized case $a=1$, on $\Omega=\R^N$ that quenches in finite time $T>0$.} Denoting the quenching point by $x_*$, they gave the following intermediate profile:
\begin{equation}\label{MZ97}
\begin{aligned}
\frac{(T-t)^{\frac{1}{r+1}}}
{1-U\!\big(x_* + z\sqrt{(T-t)|\log(T-t)|},\,t\big)}
&\sim
\left(r+1+\frac{(r+1)^2}{4r}|z|^2\right)^{-\frac{1}{r+1}}\\
&\hspace{2cm}\text{as } t\to T.
\end{aligned}
\end{equation}
Merle and Zaag also provided the final quenching profile. Namely, there exists a limiting profile $U^*$, with $U(x,t)\to U^*(x)$ as $t\to T$, such that the distance to the singular state satisfies
\begin{equation}\label{eq:MZ-final}
\begin{aligned}
1-U^*(x)
&\sim
\left(\frac{(r+1)^2}{8r}
\frac{|x-x_*|^2}{\bigl|\log |x-x_*|\bigr|}\right)^{\frac{1}{r+1}}\\
&\hspace{2cm}\text{as } x\to x_*.
\end{aligned}
\end{equation}
Later on, Huang and Zaag \cite{HZ2024} refined this construction and provided the profile of the gradient of the solution. Furthermore, in \cite{FG1993}, a general result was established for a class of solutions exhibiting the profile given in \eqref{eq:MZ-final}.

In comparison with the theoretical aspects, the numerical analysis of quenching phenomena has received relatively little attention, to the best of our knowledge. Our work is motivated by the numerical results presented in \cite{N2017} and \cite{CW2024} for blow-up solutions of heat equations, using the rescaling algorithm introduced by Berger and Kohn in \cite{BK1988}. This method is based on rescaling a portion of the solution once it reaches a certain amplitude, starting from a lower threshold (upper threshold, in quenching). In this paper, we introduce a modified version of the algorithm, in which the entire solution is rescaled, rather than just a part of it.

In this paper, we assume that $U$ quenches and focus on the normalized model with $a=1$. Moreover, we are in the one dimensional case, i.e. \revred{$N=1$ with $\Omega_0:=\Omega=(-1,1)$, where $\Omega_0$ denotes the original (unscaled) domain, and set}
$$u(x,t)=1-U(x,t).$$
Therefore, equation \eqref{EQ} becomes
\begin{equation}\label{eq}
\begin{cases}
    &\partial_t u = \partial_{x}^2 u    -
\displaystyle \frac{ 1}{  u^r }
,\ x\in(-1,1), t>0,\\
&u(1,t)=u(-1,t)=1, t>0,\\
&u(x,0)=u_0(x), x\in [-1,1],
\end{cases}
\end{equation}
where $u_0=1-U_0$. Therefore, the study of a solution when \eqref{eq:touchdown-U} occurs is equivalent to the following:
\begin{equation}\label{eq:touchdown-u}
   \underset{t\rightarrow T}{\lim } [\min_{x\in \overline\Omega}u(x,t) ]=0 \mbox{ when }t\to T.
\end{equation}

The rescaling algorithm relies on the scaling of the differential equation.  More precisely, the equation
$$
    u_t=u_{xx}-u^{-r}
$$
is invariant under
$$
    u_\lambda(x,t)=\lambda^\alpha u(\lambda^{\beta/2}x,\lambda^\beta t)
$$
provided that
\begin{equation}\label{eq:scale-relation}
\beta=-(r+1)\alpha.
\end{equation}
\revred{This statement concerns the differential equation itself, not the fixed boundary problem on $(-1,1)$ with boundary value $1$.}  On bounded intervals the same change of variables also transforms the interval and the boundary values: a stage posed on $\Omega_k$ with boundary value $\lambda^{k\alpha}$ is mapped to a stage posed on
$$
    \Omega_{k+1}=\lambda^{-\beta/2}\Omega_k
$$
with boundary value $\lambda^{(k+1)\alpha}$.  The stagewise algorithm below is formulated in this transformed-domain sense.

In this work, we investigate the convergence properties of a rescaling algorithm applied to a quenching-type problem for both equation \eqref{eq} and the associated ordinary differential equation, i.e.

\begin{equation}\label{ODE}
\begin{cases}
    &    \partial_t w =  -
\displaystyle \frac{ b}{  w^r }, t>0,\\
&w(0)=w_0,
\end{cases}
\end{equation}
where $b>0$. \revred{For the ODE case, Theorem \ref{thm:ode-upper-bound} establishes that the numerically computed quenching time converges to the exact quenching time with first-order accuracy in time.}  For the PDE, Theorem \ref{thm:convergence} proves a second-order finite-stage convergence result.  This theorem is deliberately stated for any fixed number of rescaling stages, since the exact solution is singular at the quenching time.  The full numerical quenching time is then treated in Theorem \ref{thm:pde-quenching-time-convergence} under a uniform tail comparison assumption.  This separation keeps the finite-stage consistency argument distinct from the additional tail control needed near quenching.  This paper is organized as follows:

\begin{itemize}
    \item In Section \ref{sec:ode}, we introduce the rescaling algorithm for the related ODE and establish convergence results for the numerical quenching time.
    \item In Section \ref{PDE}, we extend the rescaling algorithm to the PDE case and provide convergence results for both the numerical solution and the numerical quenching time.
    \item In Section \ref{sec:numerics}, we finally present numerical experiments illustrating the method on a simple example.
\end{itemize}
\section{On the quenching time for the ODE}\label{sec:ode}
In this section, we apply the rescaling algorithm to the ODE model \eqref{ODE} and analyze the resulting numerical quenching time.

We first introduce an Euler approximation of \eqref{ODE}. Let $\tau > 0$ be a time step and let $n \geq 0$ be a nonnegative integer. We set $t_n = n\tau$. In what follows, $W_{n}$ denotes the numerical approximation and $w_{n}=w(t_n)$ denotes the exact value. The Euler discretization of \eqref{ODE} is as follows: for $n \geq 0$,
\begin{equation}\label{eq:ode-scheme}
\begin{cases}
& \delta_t W_{n} = - \frac{b}{W_{n}^{r}},\\
&W_0=w_0,
\end{cases}
\end{equation}
where $\delta_t W_{n} $ is defined as it follows:
$$\delta_t W_{n} = \frac{W_{n+1} - W_{n}}{\tau}.$$
\begin{prp}[Monotonicity and positivity of the numerical solution]\label{prop:ode-positivity}
For a given $M\in (0,1]$, we assume that $\tau < \frac{M^{1+r}}{b}$. Consider $n\in \N$. If the following holds
$$M\leq W_{n}.$$
Then, we have
$$ 0< W_{n+1}.$$
\end{prp}
\begin{proof}
Consider $n\in \N$. From \eqref{eq:ode-scheme}, we have
    {\color{revisionred}
    \begin{equation}
        \begin{aligned}
        W_{n+1}&=W_{n}-b\tau (W_{n})^{-r}\ge M-b\tau M^{-r}>0.
        \end{aligned}
    \end{equation}
    }
This concludes the proof of Proposition \ref{prop:ode-positivity}.
\end{proof}
\subsection{Rescaling algorithm for the ordinary differential equation}
We now describe the rescaling algorithm for \eqref{ODE}. Let $\tau > 0$ be a time step and set $t_n=n\tau$.
 We introduce the function $F_{\tau}$  as a piecewise linear interpolation in temporal dimension, ensuring $F_{\tau}(t_n) = F_{n}$. The function $F_{\tau}$ can be explicitly defined for all $t \in (t_n, t_{n+1})$ as:
\begin{equation}
\begin{aligned}\label{eq:time-interp-ode}
 F_{\tau}( t) = \frac{1}{ \tau} \big[ F_{n}(t_{n+1} - t)+ F_{n+1}(t - t_n)  \big].
\end{aligned}
\end{equation}
\revred{We use the scaling relation \eqref{eq:scale-relation}, so that $\alpha=-\beta/(r+1)<0$, and set the parameters}
\begin{equation}\label{eq:ode-parameters}
\begin{aligned}
&0<M<w_0,\qquad \beta>0,\qquad M^{-\frac{1}{\alpha}}<\lambda<1,\\
&\lambda^{-\frac{\beta}{2}}\in\mathbb N,\qquad 0<\tau < \frac{M^{1+r}}{b}.
\end{aligned}
\end{equation}
We integrate the numerical solution of \eqref{eq:ode-scheme} until the first time step $\tilde n_0$ such that $W_{\tilde n_0}\le M$.  The crossing time is then located by the linear interpolant \eqref{eq:time-interp-ode}: we choose $\tilde\tau_0^*$ such that
$$
    (\tilde n_0-1)\tau\le \tilde\tau_0^*\le \tilde n_0\tau,
    \qquad
    W_\tau(\tilde\tau_0^*)=M.
$$
The next numerical stage is initialized by
$$
    W_0^{(1)}=\lambda^\alpha W_\tau(\tilde\tau_0^*)=\lambda^\alpha M,
$$
and then the same Euler recurrence \eqref{eq:ode-scheme} is applied with the same time step $\tau$.  This is equivalent to observing the original variable on a physical time scale reduced by the factor $\lambda^\beta$.  In the ODE estimates below we work directly with the numerical sequence and the weighted telescoping identity; hence no exact rescaled stage based on the numerical crossing time is needed.

The computation of $W^{(1)}$ stops at the first time step $\tilde n_1\ge0$ such that $W_{\tilde n_1}^{(1)}\le M$.  We define $\tilde\tau_1^*$ by
$$
    (\tilde n_1-1)\tau\le \tilde\tau_1^*\le \tilde n_1\tau,
    \qquad
    W_\tau^{(1)}(\tilde\tau_1^*)=M.
$$
Repeating the same procedure, the $k$-th numerical stage is stopped at the first interpolated crossing
$$
    (\tilde n_k-1)\tau\le \tilde\tau_k^*\le \tilde n_k\tau,
    \qquad
    W_\tau^{(k)}(\tilde\tau_k^*)=M,
$$
and the following stage is initialized by
$$
    W_0^{(k+1)}=\lambda^\alpha W_\tau^{(k)}(\tilde\tau_k^*)=\lambda^\alpha M.
$$
We define for all $K\in \N$, the following
\begin{equation}\label{eq:def-tilde-mu}
\tilde\mu_K=  \sum_{k=0}^K \lambda^{\beta k} \tilde\tau^*_k.
\end{equation}
Then, for all $K\in \N$, for each $ t \in  [0, \tilde\mu_K] $, the threshold is crossed only at the right endpoint of a stage:
\begin{itemize}
    \item if $0\le t<\tilde\tau^*_0$, then $W_{\tau}^{(0)}(t)> M$, while $W_{\tau}^{(0)}(\tilde\tau^*_0)=M$;
    \item if $1\le k\le K$ and $\tilde\mu_{k-1}\le t<\tilde \mu_k$, then
    $$ W^{(k)}_{\tau}\bigl( \lambda^{-\beta k}(t - \tilde \mu_{k-1})\bigr)> M;$$
    \item if $t =\tilde\mu_k$ for some $0\le k\le K$, then
    $$W^{(k)}_{\tau}\bigl(\lambda^{-\beta k}(t - \tilde\mu_{k-1})\bigr)= M,$$
    with the convention $\tilde\mu_{-1}=0$.
\end{itemize}
Therefore, the numerical solution $ W_{\tau}(t) $ of the rescaling method is defined stagewise on half-open intervals, with the endpoint values assigned by continuity:
$$
W_{\tau}(t) =
\lambda^{-\alpha k} W^{(k)}_{\tau} ( \lambda^{-\beta k}(t - \tilde\mu_{k-1}))\text{ if } t \in [\tilde\mu_{k-1},\tilde\mu_k),
$$
for $k\ge0$ and $\tilde\mu_{-1}=0$.
One can see that this solution tends to $0$ when $ K $ goes to infinity. The finite numerical quenching time is given by the following:
\begin{equation}\label{eq:tilde-T-tau}
\tilde T_{\tau} =\lim_{K \to +\infty} \tilde\mu_K= \lim_{K \to +\infty} \sum_{k=0}^K \lambda^{\beta k} \tilde\tau^*_k < +\infty.
\end{equation}

\subsection{Analysis on the quenching time for the ODE}
We now prove the convergence of the numerical quenching time.  The proof is based on a one-stage crossing estimate and on the telescoping identity induced by the scaling relation \eqref{eq:scale-relation}.

\begin{lemma}[One-stage ODE crossing estimate]\label{lem:ode-local-crossing}
\revred{Let $b,r,M>0$ be as in Proposition \ref{prop:ode-positivity}, and fix $A>M$.}  There exist constants $C_A>0$ and $\tau_A\in(0,M^{r+1}/b)$ such that the following holds.  If $M<a\le A$ and the Euler sequence
$$
    W_{n+1}=W_n-b\tau W_n^{-r},\qquad W_0=a,
$$
is stopped at the first linearly interpolated time $\widehat\tau$ for which $W_\tau(\widehat\tau)=M$, then, for $0<\tau\le\tau_A$,
\begin{equation}\label{eq:ode-local-crossing}
    \left|
    \widehat\tau-\frac{a^{r+1}-M^{r+1}}{(r+1)b}
    \right|
    \le C_A\tau .
\end{equation}
\end{lemma}
\begin{proof}
\revred{As long as $W_n>M$, one has $W_n\le a\le A$ and $W_{n+1}\le W_n-b\tau A^{-r}$. Thus the threshold is reached after finitely many steps.} Let $\widehat n$ be the first integer such that $W_{\widehat n}\le M$.  For all $0\le n\le \widehat n-1$ one has $W_n>M$.  Put $Z_n=W_n^{r+1}$.  Since
$$
    W_{n+1}=W_n\left(1-b\tau W_n^{-r-1}\right),
$$
and $W_n\ge M$, Taylor's formula on the interval $0\le x\le b\tau_A M^{-r-1}<1$ gives
$$
    Z_n-Z_{n+1}=(r+1)b\tau+\rho_n,
    \qquad |\rho_n|\le C_A\tau^2 .
$$
Summing this from $0$ to $\hat n-1$, we obtain
\begin{equation}\label{eq:Z-step}
\left|\hat n\tau-\frac{Z_0-Z_{\hat n}}{(r+1)b}\right|\leq\frac{1}{(r+1)b}\left|\sum_{n=0}^{\hat{n}-1}\rho_n\right|\le \frac{C_A\tau^2\hat n}{(r+1)b}.
\end{equation}
On the other hand, after decreasing $\tau_A$ if necessary, we may also assume
$$
    |\rho_n|\le \frac12(r+1)b\tau,
$$
for all steps before crossing.  Hence
$$
    Z_n-Z_{n+1}\ge \frac12(r+1)b\tau.
$$
Together with $0\le Z_n\le A^{r+1}$ up to the crossing step, this gives the uniform bound
$$
    \widehat n\tau\le \frac{2A^{r+1}}{(r+1)b}.
$$
With \eqref{eq:Z-step}, this yields
$$
    \left|
    \widehat n\tau-\frac{a^{r+1}-W_{\widehat n}^{r+1}}{(r+1)b}
    \right|
    \le C_A\tau .
$$
Moreover $W_{\widehat n-1}>M\ge W_{\widehat n}$. Then, we get
{\color{revisionred}
$$
    0\le M-W_{\widehat n}<W_{\widehat n-1}-W_{\widehat n}=b\tau W_{\widehat n-1}^{-r}
    \le b\tau M^{-r}.
$$
}
So
$$|W_{\widehat n}^{r+1}-M^{r+1}|=(r+1)\xi^r |W_{\widehat n}-M|\le C_A\tau.$$

Finally, the interpolated crossing time satisfies $|\widehat\tau-\widehat n\tau|\le\tau$.  Combining these estimates proves \eqref{eq:ode-local-crossing}.
\end{proof}

\begin{theorem}[\revred{Two-sided convergence of the numerical ODE quenching time}]\label{thm:ode-upper-bound}\label{thm:ode-lower-bound}
\revblue{Assume the parameter conditions \eqref{eq:ode-parameters} and the scaling relation \eqref{eq:scale-relation}. Set $A=\max\{w_0,\lambda^\alpha M\}$ and let $0<\tau\le\tau_A$, where $\tau_A$ is the threshold in Lemma \ref{lem:ode-local-crossing}, with the same $b,r,M$ and this value of $A$.}
{\color{revisionred}
Let $\tilde T=w_0^{r+1}/((r+1)b)$ be the exact quenching time of \eqref{ODE}, and let $\tilde T_\tau$ be defined by \eqref{eq:tilde-T-tau}. Then there exists $C_{\rm ode}>0$, depending only on the fixed parameters and independent of $\tau$, such that
$$
    \left|\tilde T_\tau-\tilde T\right|
    \le C_{\rm ode}\frac{\tau}{1-\lambda^\beta}.
$$
Thus, for fixed $\lambda$ and $\beta$, the numerical quenching time converges with first-order accuracy in $\tau$.
}
\end{theorem}
\begin{proof}
For the first stage set $a_0=w_0$, and for every later stage set $a_k=\lambda^\alpha M$ $(k\ge1)$. Lemma \ref{lem:ode-local-crossing} gives
$$
    \tilde\tau_k^*
    =\frac{a_k^{r+1}-M^{r+1}}{(r+1)b}
    +\varepsilon_k,
    \qquad |\varepsilon_k|\le C_{\rm ode}\tau.
$$
{\color{revisionred}
The accumulated remainder satisfies
$$
    R_\tau:=\sum_{k=0}^{\infty}\lambda^{\beta k}\varepsilon_k,
    \qquad
    |R_\tau|\le\frac{C_{\rm ode}\tau}{1-\lambda^\beta}.
$$
Therefore
\begin{align*}
\tilde T_\tau
&=\sum_{k=0}^{\infty}\lambda^{\beta k}\tilde\tau_k^*\\
&=\frac{a_0^{r+1}-M^{r+1}}{(r+1)b}+R_\tau\\
&\quad+\frac{1}{(r+1)b}
\sum_{k=1}^{\infty}\lambda^{\beta k}
\bigl((\lambda^\alpha M)^{r+1}-M^{r+1}\bigr).
\end{align*}
}
Using $\beta=-(r+1)\alpha$, we have
$$
\sum_{k=1}^{\infty}\lambda^{\beta k}
\left((\lambda^\alpha M)^{r+1}-M^{r+1}\right)
=
\sum_{k=1}^{\infty}\left(\lambda^{\beta(k-1)}-\lambda^{\beta k}\right)M^{r+1}
=M^{r+1}.
$$
Consequently the main term equals $w_0^{r+1}/((r+1)b)=\tilde T$.
\revred{Since the remainder is bounded in absolute value by $C_{\rm ode}\tau/(1-\lambda^\beta)$, both the upper and lower bounds follow.}
\end{proof}

\section{On the quenching time for the PDE}\label{PDE}
We now focus on \eqref{eq}. \revred{Let $I\ge2$ be an integer.} We discretize the domain $\Omega = (-1, 1)$ by the grid $x_i = ih$ where $-I \leq i \leq I$ and $h = 1/I$. Although the original variable in \eqref{EQ} was denoted by $U$, from this point on $U_{i,n}$ denotes the numerical approximation of the transformed variable $u(x_i,t_n)$. \revred{The notation $U_n$ stands for $(U_{-I,n}, \dots, U_{0,n}, \dots, U_{I,n})^T$.} We denote by
\begin{equation}\label{eq:fd-operators}
\begin{aligned}
&\delta_t U_{i,n} = \frac{U_{i,n+1} - U_{i,n}}{\tau},\\
&\delta_x U_{i,n} = \frac{U_{i+1,n} - U_{i-1,n}}{2h},\\
&\delta_x^2 U_{i,n} = \frac{U_{i-1,n} - 2U_{i,n} + U_{i+1,n}}{h^2}.
\end{aligned}
\end{equation}

\subsection{Discretization of the parabolic equation}
The Euler discretization of \eqref{eq} is defined as follows: for $n \geq 0$ and $-I + 1 \leq i \leq I - 1$,
\begin{align}\label{scheme}
&\delta_t U_{i,n} = \delta_x^2 U_{i,n} - \frac{1}{U_{i,n}^{r}},\\
&U_{-I,n} = U_{I,n} = 1,\\
&U_{i,0} = \phi_i,
\end{align}
where $\phi_i = u_0(x_i)$ and $\delta_t U_{i,n}$ and $\delta_x^2 U_{i,n}$ are given in \eqref{eq:fd-operators}. For all $n\in \N$, we denote by
$$
\min U_n=\min\{U_{i,n}|-I\leq i\leq I\}.
$$

\begin{prp}[Monotonicity and positivity of the numerical solution]\label{prop:pde-positivity}
For a given $M\in (0,1]$, we assume that $\tau < \min \{\frac{h^2}{2},M^{1+r}\}$
. If for all $i\in \{-I,...,I\}$, the following holds
$$
\delta_t U_{i,n}\leq 0 \text{ and }M\leq U_{i,n}.
$$
Then, for all $i\in \{-I,...,I\}$, we have
$$
\delta_t U_{i,n+1}\leq 0 \text{ and } 0< U_{i,n+1}.
$$
\end{prp}
\begin{proof}
From the boundary condition, it is clear that the result holds for $i=-I$ and $I$. Now, for $-I+1\leq i\leq I-1$, we consider $n \in \N$ and we assume that
$$
\delta_t U_{i,n}\leq 0 \text{ and }M\leq U_{i,n}, \text{ for all } i\in \{-I+1,..., I-1\}.
$$
Then, from \eqref{scheme} and the fact that $\tau < \min \{\frac{h^2}{2},M^{1+r}\}$, we have
{\color{revisionred}
\begin{equation}
\begin{aligned}
U_{i,n+1}&=(1-2\frac{\tau}{h^2})U_{i,n}+\frac{\tau}{h^2}(U_{i-1,n}+U_{i+1,n})-\tau (U_{i,n})^{-r}\\
&\ge M-\tau M^{-r}\\
&> 0.
\end{aligned}
\end{equation}
}
Now, let us show that $\delta_t U_{i,n+1}\leq 0$. For $-I+1\leq i\leq I-1$, we have that
$$
\frac{\delta_t U_{i,n+1}-\delta_t U_{i,n}}{\tau}
=\delta^2_x (\delta_t U_{i,n})
-\frac{(U_{i,n+1})^{-r}-(U_{i,n})^{-r}}{\tau}.
$$
Using Taylor expansion and the fact $U_{i,n}, U_{i,n+1}> 0$, we have that
\begin{equation}
\begin{aligned}
\frac{\delta_t U_{i,n+1}-\delta_t U_{i,n}}{\tau}
&=\delta^2_x (\delta_t U_{i,n})+r(U_{i,n})^{-r-1}\delta_t U_{i,n}\\
&\quad -\frac{r(r+1)}{2}(\xi_{i,n})^{-r-2}
\frac{(U_{i,n+1}-U_{i,n})^{2}}{\tau}\\
&\leq \delta^2_x (\delta_t U_{i,n}).
\end{aligned}
\end{equation}
for some $\xi_{i,n}$ between $U_{i,n}$ and $U_{i,n+1}$. Using the fact that $\tau < \frac{h^2}{2}$, we obtain
$$
\delta_t U_{i,n+1}
\leq \left(1-2\frac{\tau}{h^2}\right)\delta_t U_{i,n}
+\frac{\tau}{h^2}\bigl(\delta_t U_{i-1,n}+\delta_t U_{i+1,n}\bigr)
\leq 0,
$$
which concludes the proof of Proposition \ref{prop:pde-positivity}.
\end{proof}

\subsection{Rescaling algorithm for the PDE}\label{sec:pde-rescaling-algorithm}
We now describe the stagewise rescaling algorithm for \eqref{eq}. The method follows the Berger and Kohn rescaling philosophy, but it is applied to the vanishing variable $u$. In contrast with a purely piecewise-linear space-time interpolation, the present implementation uses linear interpolation only in time and a local cubic reconstruction in space at each rescaling step.

For nodal data $F_{i,n}$ we use the linear time interpolant
\begin{equation}\label{eq:time-interpolation-nodes}
 F_{i,\tau}(t)
 =\frac{1}{\tau}\bigl[F_{i,n}(t_{n+1}-t)+F_{i,n+1}(t-t_n)\bigr],
 \qquad t\in[t_n,t_{n+1}].
\end{equation}
\revred{The numerical stopping time is defined from these nodal interpolants. In general, interpolating the two endpoint minimum values gives a different function. We write $g_{k,\tau}$ and $\tau_{k,\tau}^*$ without an $h$ subscript to simplify the notation; both still depend on the spatial mesh through the numerical nodal values. At the $k$-th stage we put}
{\color{revisionred}
\begin{equation}\label{eq:nodal-interpolated-minimum}
    g_{k,\tau}(t):=\min_{-I_k\le i\le I_k}U_{i,\tau}^{(k)}(t),
\end{equation}
}
and define the numerical stopping time as the first crossing
{\color{revisionred}
\begin{equation}\label{eq:numerical-stopping-time}
    \tau_{k,\tau}^*
    :=\inf\{t>0:\ g_{k,\tau}(t)\le M\}.
\end{equation}
}
Equivalently, if $n_k$ is the first index such that $\min_iU_{i,n_k}^{(k)}\le M$, then
{\color{revisionred}
\begin{equation}\label{eq:theta-node-wise}
    \theta_k
    =\min_{\{i:\,U_{i,n_k}^{(k)}\le M\}}
    \frac{U_{i,n_k-1}^{(k)}-M}
         {U_{i,n_k-1}^{(k)}-U_{i,n_k}^{(k)}},
    \qquad
    \tau_{k,\tau}^*=(n_k-1+\theta_k)\tau.
\end{equation}
}
The minimum is taken only over nodes that cross the level $M$ in the step. The stopping profile is then
{\color{revisionred}
\begin{equation}\label{eq:W-stopping-profile}
    W_i^{(k)}:=U_{i,\tau}^{(k)}(\tau_{k,\tau}^*)
    =(1-\theta_k)U_{i,n_k-1}^{(k)}+\theta_kU_{i,n_k}^{(k)}.
\end{equation}
}
This node-wise definition guarantees that at least one grid node reaches $M$:
\begin{equation}\label{eq:min-W-equals-M}
    \min_iW_i^{(k)}=M.
\end{equation}
This property is essential in the tail comparison argument below.

Spatial values needed after a rescaling are not obtained by piecewise linear interpolation. \revred{For a grid profile $Y_i$ on $x_i=ih$, and for an interior cell $[x_i,x_{i+1}]$, let $P_i[Y]$ be the cubic Lagrange interpolation polynomial through the four values $Y_{i-1},Y_i,Y_{i+1},Y_{i+2}$.} Explicitly,
\begin{equation}\label{eq:cubic-reconstruction-polynomial}
\begin{aligned}
P_i[Y](x)
&=\frac{Y_{i+2}-3Y_{i+1}+3Y_i-Y_{i-1}}{6h^3}(x-x_i)^3\\
&\quad +\frac{Y_{i+1}-2Y_i+Y_{i-1}}{2h^2}(x-x_i)^2\\
&\quad
  -\frac{Y_{i+2}-6Y_{i+1}+3Y_i+2Y_{i-1}}{6h}(x-x_i)+Y_i.
\end{aligned}
\end{equation}
{\color{revisionred}
Equivalently, the Lagrange form of the same polynomial is
$$
P_i[Y](x)=\sum_{j\in S_i}Y_j\ell_{i,j}(x),
\qquad
\ell_{i,j}(x):=\prod_{\substack{m\in S_i\\m\ne j}}\frac{x-x_m}{x_j-x_m},
$$
where $S_i=\{i-1,i,i+1,i+2\}$ for an interior stencil. We use only the cellwise notation $P_i[Y]$ for spatial interpolation.

On a grid with endpoints $x_{-J}$ and $x_J$, we use the same
Lagrange formula on the boundary cells. At $x_{-J}$ we use the stencil
$S_{-J}:=S_{-J+1}$, and at $x_J$ and its adjacent node $x_{J-1}$ we use
$S_J:=S_{J-1}=S_{J-2}$. Explicitly,
$$
S_{-J}=\{-J,-J+1,-J+2,-J+3\},\qquad
S_{J-1}=\{J-3,J-2,J-1,J\}.
$$
At the right endpoint, $P_J[Y](x_J):=P_{J-1}[Y](x_J)=Y_J$.
The cell indices run from $-J$ to $J-1$.

These four-node stencils define the interpolation on the entire closed interval without ghost values. Because there are only three stencil shapes after scaling by $h$, their Lagrange weights satisfy a common bound
$$
\sup_{x\in[x_i,x_{i+1}]}\sum_{j\in S_i}|\ell_{i,j}(x)|\le C_L,
$$
where $C_L$ is independent of $h$, the cell, and the stage. For samples of a $C^4$ function, the Lagrange remainder is bounded by $Ch^4\|f^{(4)}\|_\infty$ on the stencil. The special difference identities used later for the rescaling defects concern interior stencils only; the finite-stage interpolation error estimates also cover the stated one-sided stencils.
}

We set
{\color{revisionred}
\begin{equation}\label{parameters}
\begin{aligned}
&M\in (0,1),\qquad \beta>0,\qquad \alpha=-\frac{\beta}{r+1},\\
&M^{-1/\alpha}<\lambda<1,\qquad q:=\lambda^{-\frac{\beta}{2}}\in\mathbb N,
\end{aligned}
\end{equation}
}
\revred{so that the integer $q$ is the mesh-alignment ratio (in particular $q\ge2$):}
\begin{equation}\label{eq:aligned-q}
    I_k=q^kI,
    \qquad \xi_{k,i}=ih,
    \qquad \eta_{k,n}=n\tau.
\end{equation}
The exact and numerical rescaled domains are
\begin{equation}\label{eq:Omega-k-definition}
    \Omega_0=(-1,1),
    \qquad
    \Omega_k=\lambda^{-\frac{\beta k}{2}}\Omega_0=(-q^k,q^k),
\end{equation}
so that the nodes $\xi_{k,i}=ih$, $-I_k\le i\le I_k$, cover $\overline\Omega_k$.
We also assume
{\color{revisionred}
$$
\delta_t U_{i,0}\leq 0,
\qquad M<U_{i,0}\leq 1,
\qquad \text{for all } -I+1\le i\le I-1,
$$
}
and
$$
    \tau < \min \left\{\frac{h^2}{2},M^{1+r}\right\}.
$$
For $k=0$ this gives $\xi_{0,i}=x_i$ and $\eta_{0,n}=t_n$. \revred{The step-size restriction is the one in Proposition \ref{prop:pde-positivity}, whose positivity estimate applies at every stage as long as its nodal values start above $M$ and its boundary values are at least $M$. Here this boundary condition is automatic, since
$U_{\pm I_k,n}^{(k)}=\lambda^{k\alpha}\ge1>M$
for $0<\lambda<1$ and $\alpha<0$. Preservation of temporal monotonicity additionally requires the nonpositive initial time differences from that proposition; it is not automatic under cubic interpolation. When later estimates need this property after rescaling, it is included explicitly among their hypotheses.}

Let $u^{(0)}=u$. {\color{revisionred}The exact stopping time $\tau_k^*$ is the first time at which the exact stage reaches $M$; hence
$$
    \min_{\xi\in\Omega_k}u^{(k)}(\xi,\tau_k^*)=M,
$$
and the next exact stage is defined by}
\begin{equation}\label{eq:exact-stage-pde}
    u^{(k+1)}(\xi_{k+1},\eta_{k+1})
    =\lambda^\alpha u^{(k)}(\lambda^{\beta/2}\xi_{k+1},\tau_k^*+\lambda^\beta\eta_{k+1}).
\end{equation}
\revred{The time shift in \eqref{eq:exact-stage-pde} is the exact time $\tau_k^*$, not the numerical time $\tau_{k,\tau}^*$. This distinction is essential for comparing the exact stages with the numerical stages.} The scale relation \eqref{eq:scale-relation} makes the differential equation invariant under this transformation. The boundary value problem is not invariant on the fixed interval $(-1,1)$; rather, the interval and boundary values are rescaled, which is why the stages are posed on $\Omega_k$ and have boundary values $\lambda^{k\alpha}$. {\color{revisionred}For the numerical initialization, set $a_i=\lambda^{\beta/2}\xi_{k+1,i}$ and choose the old cell index
$$
j(i)=\min\!\left\{\left\lfloor\frac{i}{q}\right\rfloor,\ I_k-1\right\},
\qquad -qI_k\le i\le qI_k.
$$
The next numerical initial profile is produced from the stopping profile \eqref{eq:W-stopping-profile}, evaluated at the numerical time $\tau_{k,\tau}^*$, by
\begin{equation}\label{eq:next-stage-cubic}
    U_{i,0}^{(k+1)}
    =\lambda^\alpha P_{j(i)}[W^{(k)}](a_i),
    \qquad -I_{k+1}\le i\le I_{k+1}.
\end{equation}
}
In \eqref{eq:next-stage-cubic}, the new node $\xi_{k+1,i}$ corresponds to the old spatial point $a_i=\xi_{k+1,i}/q$. We evaluate the cubic interpolant of the old stopping profile at this point and multiply by $\lambda^\alpha$. In particular, at a retained old node $i=qj$, the interpolation property gives
$$
U_{qj,0}^{(k+1)}=\lambda^\alpha W_j^{(k)}.
$$

The boundary values in the rescaled variables are
$$
    U_{I_{k+1},n}^{(k+1)}=U_{-I_{k+1},n}^{(k+1)}=\lambda^{(k+1)\alpha}.
$$

The stopping condition is imposed on nodal values. \revred{Each numerical stage is required to start with $\min_iU_{i,0}^{(k)}>M$. For any fixed number of stages this follows, for sufficiently small $h$, from the finite-stage comparison proved below, since the next exact stage has minimum $\lambda^\alpha M>M$. For an infinite sequence of stages, admissibility of the restarts is part of the tail hypotheses. The positivity condition needed for the cubic reconstruction is also stated explicitly below.}

We define the numerical accumulated stopping times by
{\color{revisionred}
\begin{equation}\label{eq:def-mu-h-tau}
    \mu_{K,h,\tau}:=\sum_{k=0}^K\lambda^{\beta k}\tau_{k,\tau}^*,
\end{equation}
}
and the exact accumulated stopping times by
\begin{equation}\label{eq:exact-muK}
    \mu_K:=\sum_{k=0}^K\lambda^{\beta k}\tau_k^*.
\end{equation}
{\color{revisionred}
For a point $\xi$ in the cell $[\xi_{k,j},\xi_{k,j+1}]$, define the spatially reconstructed stage value by
$$
U_{h,\tau}^{(k)}(\xi,t):=P_j[U_{\cdot,\tau}^{(k)}(t)](\xi).
$$
}
For $\mu_{-1,h,\tau}=0$, the physical numerical solution is defined stagewise by
{\color{revisionred}
$$
\begin{aligned}
U_{h,\tau}(x,t)
&=\lambda^{-\alpha k}U_{h,\tau}^{(k)}
\left(\lambda^{-\frac{\beta}{2}k}x,\lambda^{-\beta k}(t-\mu_{k-1,h,\tau})\right),\\
&\hspace{3cm}t\in[\mu_{k-1,h,\tau},\mu_{k,h,\tau}).
\end{aligned}
$$
The numerical quenching time, provided that all stages are admissible, is
\begin{equation}\label{T_tau,h}
    T_{h,\tau}:=\lim_{K\to\infty}\mu_{K,h,\tau}
    =\sum_{k=0}^{\infty}\lambda^{\beta k}\tau_{k,\tau}^*.
\end{equation}
}
We now record the discrete differential inequality needed in the sequel. The point is that the rescaling step does not preserve the inequality with no error; the temporal interpolation error and the cubic spatial interpolation error have to be accumulated with the correct scaling factor.

\subsection{Cubic spatial interpolation and accumulated defects}
{\color{revisionred}
For the local interpolation estimates, we retain the fixed core
\begin{equation}\label{def:Omega-lambda}
    \Omega_\lambda:=\big[-\lambda^{2\beta},\lambda^{2\beta}\big]\subset\Omega_0.
\end{equation}
Here the original half-length is $1$. This core is distinct from the expanding stage domain $\Omega_k$. Every local interpolation stencil below must be contained in $\Omega_\lambda$. The exponent $\beta>0$ remains arbitrary; no specialization to $\beta=2$ is made in this analysis.

We use the integer mesh ratio $q$ already defined in \eqref{parameters}, and introduce
\begin{equation}\label{def:q-mu-lambda}
    \mu_\lambda:=\lambda^{\alpha+\beta}
    =\lambda^{\frac{\beta r}{1+r}}<1.
\end{equation}
The relation $\beta=-(r+1)\alpha$ implies $\lambda^{\alpha+\beta}=\lambda^{-\alpha r}$. Thus $\mu_\lambda$ is the factor by which a spatial rescaling multiplies an additive defect. The symbol $q$ is needed only for the alignment of the two grids and for stencil coefficients.

We also record explicitly the constants used for evolution inside a stage. Assume that all nodal values used up to and including the crossing step are at least $m_\lambda>0$, independently of $k,h,\tau$. For example, if a stage starts above $M$, the positivity argument in Proposition \ref{prop:pde-positivity}, with $\tau\le M^{r+1}/2$ and $\tau/h^2\le1/2$, gives the lower bound $m_\lambda=M/2$ at its first crossing step. Assume in addition that the stage lengths satisfy
\begin{equation}\label{ass:uniform-stage-length}
    n_k\tau\le T_\lambda
\end{equation}
uniformly in $k$, and define
\begin{equation}\label{def:A-lambda}
    L_\lambda:=rm_\lambda^{-r-1},
    \qquad
    A_\lambda:=\exp(L_\lambda T_\lambda).
\end{equation}
The constant $A_\lambda$ bounds the amplification of a positive additive defect during the Euler evolution; see Lemma \ref{lem:defect-propagation-stage}. This factor is needed also when the final estimate is evaluated only at the stopping time.

For a grid profile $Y_i$, define
}
\begin{equation}\label{def:fourth-difference}
\delta_x^4Y_i:=\frac{Y_{i+2}-4Y_{i+1}+6Y_i-4Y_{i-1}+Y_{i-2}}{h^4}.
\end{equation}
\revred{The cubic Lagrange polynomials in \eqref{eq:cubic-reconstruction-polynomial} have the following properties.}
\begin{lemma}[Local identities for the cubic reconstruction]\label{lemma:cubic-identities}
{\color{revisionred}
Let $P_i[Y]$ be defined by \eqref{eq:cubic-reconstruction-polynomial}, and let $\Delta\ne0$. Then
\begin{align}
P_i[Y](x_j)&=Y_j,
\qquad j=i-1,i,i+1,i+2,\label{eq:cubic-identity-interpolation}\\
\frac{P_i[Y](x+\Delta)-2P_i[Y](x)+P_i[Y](x-\Delta)}{\Delta^2}
&=P_i[Y]''(x),\label{eq:cubic-identity-second}\\
P_i[Y]''(x_i+\theta h)
&=(1-\theta)\delta_x^2Y_i+\theta\delta_x^2Y_{i+1},
\label{eq:cubic-identity-linear-second}\\
&\hspace{2cm}0\le\theta\le1.\notag
\end{align}
}
and
\begin{equation}\label{eq:cubic-identity-jump}
P_i[Y](x)-P_{i-1}[Y](x)
=\frac{h\,\delta_x^4Y_i}{6}(x-x_{i-1})(x-x_i)(x-x_{i+1}).
\end{equation}
\end{lemma}

\begin{proof}
{\color{revisionred}
The first identity is the defining interpolation property. The Taylor expansion of a cubic polynomial is exact through degree three; the odd terms cancel in its centered second difference, proving \eqref{eq:cubic-identity-second}. At $x_i$ and $x_{i+1}$ this identity with $\Delta=h$ gives respectively $\delta_x^2Y_i$ and $\delta_x^2Y_{i+1}$. Since $P_i[Y]''$ is affine, its value between these endpoints is \eqref{eq:cubic-identity-linear-second}.

Finally, $P_i[Y]-P_{i-1}[Y]$ vanishes at the three common interpolation nodes $x_{i-1},x_i,x_{i+1}$. It is therefore a constant multiple of their three linear factors. Its leading coefficient is
$$
\frac{Y_{i+2}-4Y_{i+1}+6Y_i-4Y_{i-1}+Y_{i-2}}{6h^3}
=\frac h6\delta_x^4Y_i,
$$
which proves \eqref{eq:cubic-identity-jump}.
}
\end{proof}

\begin{lemma}[Temporal interpolation defect at a stopping time]\label{lemma:temporal-defect}
\revblue{Fix $c\in[0,1]$ and assume the uniform positive lower bound $m_\lambda$ used in \eqref{def:A-lambda}.}
{\color{revisionred}
Let $\widetilde{\mathcal E}_k\ge0$, and suppose that, for $n=n_k-1,n_k$ and all required interior nodes,
\begin{equation}\label{ineq:previous-stage-with-Ek}
    \delta_x^2U_{i,n}^{(k)}
    \le
    (1-c)\left(U_{i,n}^{(k)}\right)^{-r}
    +\widetilde{\mathcal E}_k.
\end{equation}
Let $W_i^{(k)}=U_{i,\tau}^{(k)}(\tau_{k,\tau}^*)$ be the stopping profile in \eqref{eq:W-stopping-profile}. For $n_k\ge1$, set
\begin{equation}\label{def:Theta-k}
B_k:=\max_i\left|\delta_tU_{i,n_k-1}^{(k)}\right|,
\qquad
\Theta_k:=\frac{(1-c)r(r+1)}{8}m_\lambda^{-r-2}B_k^2.
\end{equation}
If no time step is taken, set $\Theta_k=0$.}
Then
{\color{revisionred}
\begin{equation}\label{ineq:nodal-time-defect}
    \delta_x^2W_i^{(k)}
    \le
    (1-c)\left(W_i^{(k)}\right)^{-r}
    +\widetilde{\mathcal E}_k+\Theta_k\tau^2.
\end{equation}
}
\end{lemma}

\begin{proof}
{\color{revisionred}
Write $\theta=\theta_k$ from \eqref{eq:theta-node-wise}, and put $a=U_{i,n_k-1}^{(k)}$, $b=U_{i,n_k}^{(k)}$. Since the same $\theta$ is used at every node, taking $\delta_x^2$ commutes with this time interpolation. Thus \eqref{ineq:previous-stage-with-Ek} gives
$$
\delta_x^2W_i^{(k)}
\le
(1-c)\bigl((1-\theta)a^{-r}+\theta b^{-r}\bigr)
+\widetilde{\mathcal E}_k.
$$
For $f(s)=s^{-r}$, Taylor's interpolation remainder, together with $a,b\ge m_\lambda$, yields
$$
0\le (1-\theta)f(a)+\theta f(b)-f((1-\theta)a+\theta b)
\le\frac{r(r+1)}{8}m_\lambda^{-r-2}(b-a)^2.
$$
Here the factor $1/8$ comes from $\theta(1-\theta)/2\le1/8$. Since $b-a=\tau\delta_tU_{i,n_k-1}^{(k)}$, the additional term is at most $\Theta_k\tau^2$.
}
\end{proof}

\begin{lemma}[Cubic rescaling step and spatial defect]\label{lem:cubic-rescaling-defect}
\revblue{Use the same $c\in[0,1]$ as in Lemma \ref{lemma:temporal-defect}, and the aligned grids of \eqref{eq:aligned-q}.}
{\color{revisionred}
Suppose that
\begin{equation}\label{ineq:nodal-profile-with-defect}
    \delta_x^2W_i^{(k)}
    \le
    (1-c)\left(W_i^{(k)}\right)^{-r}
    +\widehat{\mathcal E}_k
\end{equation}
at all old-grid nodes needed by the interpolation stencils, where $\widehat{\mathcal E}_k\ge0$. Assume that the cubic polynomials on these cells satisfy
\begin{equation}\label{ass:local-cubic-bounds}
    P_i[W^{(k)}]\ge M_{0,k}>0,
    \qquad
    \|(P_i[W^{(k)}])'\|_\infty\le C_{1,k},
    \qquad
    \|(P_i[W^{(k)}])''\|_\infty\le C_{2,k}.
\end{equation}
}
Define
{\color{revisionred}
\begin{equation}\label{def:K1-lambda}
    K_{1,k}^{\lambda}
    :=\max_i\left|\delta_x^4W_i^{(k)}\right|,
\end{equation}
where the maximum includes all required five-point stencils contained in $\Omega_\lambda$. For
$$
F_{i,k}(\theta):=\bigl(P_i[W^{(k)}](\xi_{k,i}+\theta h)\bigr)^{-r},
$$
choose $C_{F,k}^\lambda$ such that
\begin{equation}\label{def:CF-lambda}
    |F_{i,k}''(\theta)|\le C_{F,k}^\lambda h^2,
    \qquad0\le\theta\le1.
\end{equation}
These derivatives are with respect to $\theta$. The chain rule shows that \eqref{ass:local-cubic-bounds} permits the choice
}
\begin{equation}\label{def:CF-example}
C_{F,k}^{\lambda}
=r(r+1)M_{0,k}^{-r-2}C_{1,k}^2
+rM_{0,k}^{-r-1}C_{2,k}.
\end{equation}
Set
{\color{revisionred}
\begin{equation}\label{def:Dk-lambda}
    D_k:=\mu_\lambda\max\left\{
    \frac{q^2-1}{6q}K_{1,k}^{\lambda},
    \frac{1-c}{8}C_{F,k}^{\lambda}
    \right\}.
\end{equation}
Then every next-stage node whose required old-grid stencils satisfy the preceding hypotheses obeys
\begin{equation}\label{ineq:next-stage-with-Ek}
    \delta_x^2U_{j,0}^{(k+1)}
    \le
    (1-c)\left(U_{j,0}^{(k+1)}\right)^{-r}
    +\mu_\lambda\widehat{\mathcal E}_k+D_kh^2.
\end{equation}
}
\end{lemma}

\begin{proof}
{\color{revisionred}
Write $j=iq+l$, where $0\le l\le q-1$, and put $\theta=l/q$. The old spatial point corresponding to $\xi_{k+1,j}$ is $\xi_{k,i}+\theta h$. A new-grid increment $h$ corresponds to an old-grid increment $h/q$. Consequently the multiplier in a second difference is $\lambda^\alpha/q^2=\lambda^{\alpha+\beta}=\mu_\lambda$. Moreover,
$$
\mu_\lambda\bigl(P_i[W^{(k)}](\xi)\bigr)^{-r}
=\bigl(\lambda^\alpha P_i[W^{(k)}](\xi)\bigr)^{-r}.
$$

For $1\le l\le q-1$, the three values in the new second difference belong to a single old interpolation cell, including its endpoints. Lemma \ref{lemma:cubic-identities} therefore gives
$$
\delta_x^2U_{iq+l,0}^{(k+1)}
=\mu_\lambda\bigl((1-\theta)\delta_x^2W_i^{(k)}
+\theta\delta_x^2W_{i+1}^{(k)}\bigr).
$$
Using \eqref{ineq:nodal-profile-with-defect} and then linear interpolation of $F_{i,k}$, we obtain
$$
\begin{aligned}
\delta_x^2U_{iq+l,0}^{(k+1)}
&\le(1-c)\mu_\lambda\bigl((1-\theta)F_{i,k}(0)+\theta F_{i,k}(1)\bigr)
+\mu_\lambda\widehat{\mathcal E}_k\\
&\le(1-c)\left(U_{iq+l,0}^{(k+1)}\right)^{-r}
+\mu_\lambda\widehat{\mathcal E}_k
+\frac{1-c}{8}\mu_\lambda C_{F,k}^{\lambda}h^2.
\end{aligned}
$$
Indeed, the difference between the endpoint linear interpolant and $F_{i,k}(\theta)$ equals $\theta(1-\theta)F_{i,k}''(\eta)/2$ for some $\eta\in(0,1)$.

For $l=0$, the new second difference straddles two old interpolation cells.
Put $x=\xi_{k,i}$ and, for this calculation, abbreviate $P_j=P_j[W^{(k)}]$. The three new nodal values are
$$
\begin{aligned}
U_{iq+1,0}^{(k+1)}&=\lambda^\alpha P_i(x+h/q),\\
U_{iq,0}^{(k+1)}&=\lambda^\alpha P_i(x),\\
U_{iq-1,0}^{(k+1)}&=\lambda^\alpha P_{i-1}(x-h/q).
\end{aligned}
$$
Substitute these values into the centered second difference and add and subtract $\lambda^\alpha P_i(x-h/q)/h^2$. By \eqref{eq:cubic-identity-second} and \eqref{eq:cubic-identity-linear-second}, the single-polynomial part is
$$
\begin{aligned}
&\frac{\lambda^\alpha}{h^2}
\bigl[P_i(x+h/q)-2P_i(x)+P_i(x-h/q)\bigr]\\
&\qquad=\frac{\lambda^\alpha}{q^2}P_i''(x)
=\mu_\lambda\delta_x^2W_i^{(k)},
\end{aligned}
$$
because $q^{-2}=\lambda^\beta$. The remaining term replaces the left value of $P_i$ by the actual value of $P_{i-1}$.

Its correction to the single-polynomial formula is
$$
\begin{aligned}
\delta_x^2U_{iq,0}^{(k+1)}
&=\mu_\lambda\delta_x^2W_i^{(k)}\\
&\quad+\frac{\lambda^\alpha}{h^2}
\left[
P_{i-1}[W^{(k)}]\left(\xi_{k,i}-\frac hq\right)
-P_i[W^{(k)}]\left(\xi_{k,i}-\frac hq\right)
\right].
\end{aligned}
$$
}
{\color{revisionblue}
To evaluate the correction, put $x=\xi_{k,i}-h/q$ in \eqref{eq:cubic-identity-jump}. The product of its three factors is
$$
\begin{aligned}
&(x-\xi_{k,i-1})(x-\xi_{k,i})(x-\xi_{k,i+1})\\
&\qquad=\frac{q-1}{q}h
\left(-\frac hq\right)
\left(-\frac{q+1}{q}h\right)
=\frac{q^2-1}{q^3}h^3.
\end{aligned}
$$
The identity \eqref{eq:cubic-identity-jump} is written for $P_i-P_{i-1}$; reversing the order of the two polynomials therefore gives
$$
P_{i-1}[W^{(k)}]\left(\xi_{k,i}-\frac hq\right)
-P_i[W^{(k)}]\left(\xi_{k,i}-\frac hq\right)
=-\frac{q^2-1}{6q^3}h^4\delta_x^4W_i^{(k)}.
$$
}
\revred{Its contribution in absolute value is at most $\mu_\lambda(q^2-1)K_{1,k}^\lambda h^2/(6q)$. Combining the two cases proves \eqref{ineq:next-stage-with-Ek}.}
\end{proof}

\begin{lemma}[Propagation of an additive defect inside one stage]\label{lem:defect-propagation-stage}
\revred{Fix a stage $k$ and $c\in[0,1]$, and let $L_\lambda,A_\lambda$ be defined by \eqref{def:A-lambda}.}
\revblue{Assume the nodal lower bound $m_\lambda$ used in \eqref{def:A-lambda}, the stage-length bound \eqref{ass:uniform-stage-length}, and the CFL bound $\tau/h^2\le1/2$.}
{\color{revisionred}
Consider a consecutive set of interior grid indices $\{a,\ldots,b\}$ on which all required stencils are defined. Suppose that
\begin{equation}\label{ineq:defect-at-initial-time}
    \delta_x^2U_{i,0}^{(k)}
    \le
    (1-c)\left(U_{i,0}^{(k)}\right)^{-r}
    +\mathcal E_k,
    \qquad a\le i\le b,
\end{equation}
where $\mathcal E_k\ge0$. Define
\begin{equation}\label{def J}
    J_{i,n}^{(k)}
    :=\delta_x^2U_{i,n}^{(k)}
    -(1-c)\left(U_{i,n}^{(k)}\right)^{-r}.
\end{equation}
At the two lateral endpoints of this local index interval, assume additionally that
\begin{equation}\label{eq:lateral-defect-bound}
    J_{a,n}^{(k)},J_{b,n}^{(k)}
    \le(1+\tau L_\lambda)^n\mathcal E_k,
    \qquad0\le n\le n_k.
\end{equation}
Then
\begin{equation}\label{ineq:defect-propagated}
    \delta_x^2U_{i,n}^{(k)}
    \le
    (1-c)\left(U_{i,n}^{(k)}\right)^{-r}
    +A_\lambda\mathcal E_k,
    \quad a\le i\le b,\quad0\le n\le n_k.
\end{equation}
}
\end{lemma}

\begin{proof}
{\color{revisionred}
The convexity calculation in Appendix \ref{proof max prin J} gives, for $a<i<b$,
\begin{equation}\label{MaxPrin J}
    \delta_tJ_{i,n}^{(k)}-\delta_x^2J_{i,n}^{(k)}
    \le
    r\left(U_{i,n}^{(k)}\right)^{-r-1}J_{i,n}^{(k)}.
\end{equation}
Put $\nu=\tau/h^2$ and $a_{i,n}=r(U_{i,n}^{(k)})^{-r-1}$, so that $0\le a_{i,n}\le L_\lambda$. Rewriting \eqref{MaxPrin J} as an explicit update gives
$$
J_{i,n+1}^{(k)}
\le(1-2\nu+\tau a_{i,n})J_{i,n}^{(k)}
+\nu J_{i-1,n}^{(k)}+\nu J_{i+1,n}^{(k)}.
$$
All three coefficients are nonnegative. Let $B_n=(1+\tau L_\lambda)^n\mathcal E_k\ge0$. If $J_{i,n}^{(k)}\le B_n$ at every node, the preceding update yields
$$
J_{i,n+1}^{(k)}
\le(1+\tau a_{i,n})B_n
\le(1+\tau L_\lambda)B_n=B_{n+1}
$$
at every interior node. The lateral nodes satisfy the same bound by \eqref{eq:lateral-defect-bound}. Induction, starting from \eqref{ineq:defect-at-initial-time}, proves
$$
\max_{a\le i\le b}J_{i,n}^{(k)}
\le(1+\tau L_\lambda)^n\mathcal E_k
\le e^{L_\lambda n\tau}\mathcal E_k
\le A_\lambda\mathcal E_k.
$$
}
This proves \eqref{ineq:defect-propagated}.
\revred{The middle inequality follows from $1+x\le e^x$ for $x\ge0$. This argument uses a nonnegative upper barrier; it does not assume that $\max_iJ_{i,n}^{(k)}$ is nonnegative.}
\end{proof}

\revred{The lateral hypothesis \eqref{eq:lateral-defect-bound} is essential for a local maximum principle. Restricting a differential inequality to a core interval does not by itself control values entering from its two endpoints.}

\begin{lemma}[Propagation with no additive defect]\label{lem:no-additive-defect}
\revblue{Assume the same positivity, CFL, stage-length, and lateral comparison hypotheses as in Lemma \ref{lem:defect-propagation-stage}, with $c$ replaced by $\bar c\in[0,1]$ and $\mathcal E_k=0$.}
{\color{revisionred}
For the corresponding $J$, the lateral condition is
$$
J_{a,n}^{(k)},J_{b,n}^{(k)}\le0.
$$
If
$$
\delta_x^2U_{i,0}^{(k)}
\le(1-\bar c)(U_{i,0}^{(k)})^{-r},
\qquad a\le i\le b,
$$
then
$$
\delta_x^2U_{i,n}^{(k)}
\le(1-\bar c)(U_{i,n}^{(k)})^{-r},
\qquad a\le i\le b,\qquad0\le n\le n_k.
$$
}
\end{lemma}

\begin{proof}
\revred{Apply the nonnegative-barrier proof of Lemma \ref{lem:defect-propagation-stage} with the zero barrier.}
\end{proof}

\begin{prp}[\revred{Conditional accumulation of rescaling defects}]\label{prop:accumulated-defect}
\revblue{Fix $c\in(0,1]$. At every stage assume the hypotheses of Lemmas \ref{lemma:temporal-defect}, \ref{lem:cubic-rescaling-defect}, and \ref{lem:defect-propagation-stage}, with the same $m_\lambda,T_\lambda$, and CFL bound. The controlled local index intervals are required to be compatible with the rescaling: every old-grid stencil used to initialize a controlled next-stage node must lie among the previously controlled nodes. Impose the lateral bounds \eqref{eq:lateral-defect-bound} at every stage, with the defects defined below.}
{\color{revisionred}
Suppose initially that
\begin{equation}\label{eq:maximum-principle-assumption}
    \delta_x^2U_{i,0}^{(0)}
    \le(1-c)(U_{i,0}^{(0)})^{-r}
\end{equation}
at all controlled nodes. Let $D_k$ and $\Theta_k$ be given by \eqref{def:Dk-lambda} and \eqref{def:Theta-k}, and set
\begin{equation}\label{def:accumulated-defect}
    \mathcal E_0:=0,
    \qquad\rho_\lambda:=\mu_\lambda A_\lambda,
    \qquad\mathcal E_{k+1}:=\rho_\lambda\mathcal E_k+\mu_\lambda\Theta_k\tau^2+D_kh^2.
\end{equation}
Then
\begin{equation}\label{ineq:stage-with-Ek}
    \delta_x^2U_{i,0}^{(k)}
    \le(1-c)(U_{i,0}^{(k)})^{-r}+\mathcal E_k
\end{equation}
and
\begin{equation}\label{ineq:stage-time-with-Ek}
    \delta_x^2U_{i,n}^{(k)}
    \le(1-c)(U_{i,n}^{(k)})^{-r}+A_\lambda\mathcal E_k,
    \qquad0\le n\le n_k,
\end{equation}
on these controlled nodes.}
Moreover,
{\color{revisionred}
\begin{equation}\label{eq:Ek-explicit}
    \mathcal E_k
    =\sum_{j=0}^{k-1}\rho_\lambda^{k-1-j}
      \bigl(\mu_\lambda\Theta_j\tau^2+D_jh^2\bigr).
\end{equation}
If, in addition,
\begin{equation}\label{ass:rho-less-one}
    \rho_\lambda<1,
\end{equation}
and $D_k\le D_\lambda$, $\Theta_k\le\Theta_\lambda$ uniformly in $k$, then
\begin{equation}\label{eq:Ek-uniform}
    \mathcal E_k
    \le\frac{D_\lambda h^2+\mu_\lambda\Theta_\lambda\tau^2}{1-\rho_\lambda},
    \qquad k\ge0.
\end{equation}
If also $U_{i,n}^{(k)}\le C_\lambda$ at the controlled nodes and
\begin{equation}\label{eq:smallness-accumulated-defect}
    \frac{A_\lambda(D_\lambda h^2+\mu_\lambda\Theta_\lambda\tau^2)}{1-\rho_\lambda}
    \le\frac c2C_\lambda^{-r},
\end{equation}
then $c'=c/2$ satisfies
\begin{equation}\label{ineq:absorbed-cprime}
    \delta_x^2U_{i,n}^{(k)}
    \le(1-c')(U_{i,n}^{(k)})^{-r}
\end{equation}
throughout every controlled stage.
}
\end{prp}

\begin{proof}
{\color{revisionred}
Assume \eqref{ineq:stage-with-Ek} at stage $k$. Lemma \ref{lem:defect-propagation-stage} gives the endpoint defect $A_\lambda\mathcal E_k$ at the two time levels used for stopping. Lemma \ref{lemma:temporal-defect} increases this to $\widehat{\mathcal E}_k=A_\lambda\mathcal E_k+\Theta_k\tau^2$. Lemma \ref{lem:cubic-rescaling-defect} then gives at the next initialization the defect
$$
\mu_\lambda\widehat{\mathcal E}_k+D_kh^2
=\rho_\lambda\mathcal E_k+\mu_\lambda\Theta_k\tau^2+D_kh^2
=\mathcal E_{k+1}.
$$
This proves both stage estimates by induction. Unrolling this recurrence proves \eqref{eq:Ek-explicit}. If $\rho_\lambda<1$, the sum of its weights is at most $1/(1-\rho_\lambda)$, proving \eqref{eq:Ek-uniform}. Finally, $U_{i,n}^{(k)}\le C_\lambda$ implies $(U_{i,n}^{(k)})^{-r}\ge C_\lambda^{-r}$. Thus \eqref{eq:smallness-accumulated-defect} bounds the propagated defect by $(c/2)(U_{i,n}^{(k)})^{-r}$, which can be absorbed into the coefficient $1-c$ in \eqref{ineq:stage-time-with-Ek}.
}
\end{proof}

\begin{remark}[Raw sums versus accumulated defects]\label{rem:raw-sums}
\revred{Even when $D_k\le D_\lambda$, the unweighted sum $\sum_{k\ge0}D_kh^2$ need not converge for fixed $h>0$. The quantity in Proposition \ref{prop:accumulated-defect} is instead the weighted sum \eqref{eq:Ek-explicit}: a defect created at stage $j$ contributes at stage $k$ with weight $\rho_\lambda^{k-1-j}$. The factor $\mu_\lambda$ comes from rescaling and $A_\lambda$ from evolution between rescalings. Therefore uniform bounds on $D_k$ and $\Theta_k$ give a uniform defect of order $h^2+\tau^2$ only under the additional contraction condition $\rho_\lambda<1$. Uniform boundedness alone does not prove this conclusion.}
\end{remark}

\begin{remark}[\revblue{Uniform control of the interpolation constants}]\label{rem:uniform-cubic-constants}
{\color{revisionblue}
The uniform estimate in Proposition \ref{prop:accumulated-defect} uses bounds on the interpolation constants at every stage. The constant $K_{1,k}^\lambda$ measures the fourth spatial difference of the stopping profile, and $C_{F,k}^\lambda$ measures the curvature of the reaction evaluated on its cubic interpolants. These quantities require control during both the Euler evolution and the interpolation step.

For the interpolation step, direct substitution into the five-point fourth difference gives
$$
K_{1,k+1}^{\mathrm{init}}\le\vartheta_qK_{1,k}^{\mathrm{end}},
\qquad
\vartheta_q:=\frac{q^2+2}{3}q^{\frac{2}{1+r}-3}.
$$
Here the superscripts distinguish the values just before and just after rescaling. The factor is less than $1$ when $q=2$, and also when $r\ge1$ and $q\ge2$. For the admissible values $r=0.1$ and $q=4$, it is approximately $1.1658$. Thus the size of this interpolation factor should be checked for the chosen parameters, and the intervening Euler evolution requires its own estimate.

To bound $K_{1,k}^\lambda$ uniformly, one can bound the fourth derivatives of the exact stopping profiles and the fourth differences of the numerical errors uniformly. The latter is a discrete derivative estimate. For a grid error $e$,
$$
\|\delta_x^4e\|_\infty\le16h^{-4}\|e\|_\infty,
$$
so an $O(h^2)$ maximum-norm error permits a fourth-difference bound of order $h^{-2}$. For $C_{F,k}^\lambda$, uniform positivity and uniform first- and second-derivative bounds for the cubics give a uniform bound directly from \eqref{def:CF-example}. Establishing these estimates uniformly over all stages is an additional regularity requirement; the finite-stage theorem below provides estimates for each fixed number of stages.
}
\end{remark}

\begin{remark}[\revblue{A conditional route to decay of fourth differences}]\label{rem:decaying-fourth-difference}
{\color{revisionblue}
A decay estimate $K_{1,k}^\lambda=O((k+1)^{-2})$ would make the coarse-node contribution to the sum $\sum_kD_kh^2$ finite for fixed $h$. The following calculation explains what is needed to obtain such decay.

Suppose, as an additional hypothesis, that the exact stopping profile admits
$$
u^{(k)}(\xi,\tau_k^*)
=M+a\frac{\xi^2}{L_k}+b\frac{\xi^4}{L_k^2}+R_k(\xi),
\qquad
L_k=|\log(T-\mu_k)|,
$$
with $\|R_k\|_{C^4(\Omega_\lambda)}\le CL_k^{-3}$. The fourth difference annihilates the constant and quadratic terms and maps $\xi^4$ to $24$. Consequently the exact sampled profile has fourth difference $24b/L_k^2+O(L_k^{-3})$. If the numerical stopping profile differs from this sampled profile by an error whose fourth difference is $O(L_k^{-2})$, then $K_{1,k}^\lambda=O(L_k^{-2})$. The additional relation $L_k\asymp k+1$ gives the stated decay in the stage number. This argument uses the displayed derivative bound on the remainder and the fourth-difference bound on the numerical error.

For the full spatial-defect sum, the contribution $(1-c)C_{F,k}^\lambda$ must also be summable. The quadratic term in the displayed smooth profile can contribute a term of order $L_k^{-1}$ to the bound in \eqref{def:Dk-lambda}, so this part requires a further decay estimate when $c<1$; its coefficient vanishes when $c=1$. Proposition \ref{prop:accumulated-defect} provides a separate route to a uniform accumulated defect through uniform source bounds and the contraction condition $\rho_\lambda<1$.
}
\end{remark}

\begin{remark}[\revblue{Checking the contraction condition}]\label{rem:rho-condition}
{\color{revisionblue}
The uniform geometric estimate in Proposition \ref{prop:accumulated-defect} requires its local lateral bounds, uniform interpolation bounds, and the condition $\rho_\lambda=\mu_\lambda A_\lambda<1$. The factor $\mu_\lambda$ describes rescaling, while $A_\lambda$ bounds growth during the preceding Euler evolution. Both factors enter the contraction test.

A scalar cycle shows the size of the present amplification bound. For $w'=-w^{-r}$, the exact time to decrease from $\lambda^\alpha M$ to $M$ is
$$
t_*=\frac{M^{r+1}}{r+1}(\lambda^{-\beta}-1).
$$
Using $m_\lambda\le M$ and $T_\lambda\ge t_*$ in the definition of $A_\lambda$, and writing $z=\lambda^{-\beta}>1$, gives
$$
\log(\mu_\lambda A_\lambda)
\ge\frac r{r+1}(z-1-\log z)>0.
$$
Thus this amplification bound gives $\rho_\lambda>1$ for the scalar cycle. A contraction proof through this approach requires a sharper estimate for evolution within a stage.

The finite-stage convergence theorem below uses Assumptions \ref{ass:finite-stage} and \ref{ass:cfl-reconstruction}. Convergence of the full numerical quenching time additionally requires the uniform tail comparison in Assumption \ref{ass:tail-comparison}. Proposition \ref{prop:accumulated-defect} supplies a sufficient route to that comparison when all of its hypotheses, including contraction and defect absorption, are verified on a region containing the relevant nodal minima.
}
\end{remark}

\subsection{Finite-stage convergence}\label{subsec:conv}
\revred{The exact solution is singular at the quenching time, so the result in this subsection concerns a fixed finite number of stages. Constants may depend on that number, but not on $h$ or $\tau$. The exact first crossing time is denoted by $\tau_k^*$, and the numerical first crossing time by $\tau_{k,\tau}^*$, as in \eqref{eq:numerical-stopping-time}. The spatial mesh dependence of the latter notation is suppressed, not removed.}

\begin{assumption}[Finite-stage regularity and transversality]\label{ass:finite-stage}
Fix $K\in\mathbb N$.
\revred{There exists $\delta_K>0$, with $\delta_K<\min_{0\le k\le K}\tau_k^*$, such that, for every $0\le k\le K$, the exact stage $u^{(k)}$ is defined on
$$
    Q_{k,K}:=\overline\Omega_k\times[0,\tau_k^*+2\delta_K]
$$
and is bounded in $C^{4,2}(Q_{k,K})$. There are constants $0<m_K<B_K<\infty$, independent of $h,\tau$, such that
$$
    m_K\le u^{(k)}\le B_K
    \qquad\text{on }Q_{k,K}.
$$
In particular, its spatial derivatives through order four and its time derivatives through order two are bounded on these sets.}
Put
\begin{equation}\label{eq:exact-min-gk}
    g_k(\eta):=\min_{\xi\in\overline\Omega_k}u^{(k)}(\xi,\eta).
\end{equation}
\revred{The time $\tau_k^*>0$ is the first crossing of $M$: $g_k(t)>M$ for $0\le t<\tau_k^*$ and $g_k(\tau_k^*)=M$. Moreover, there exists $\sigma_K>0$ such that, for $0\le s\le\delta_K$,}
\begin{equation}\label{eq:transversality}
    g_k(\tau_k^*-s)\ge M+\sigma_Ks,
    \qquad
    g_k(\tau_k^*+s)\le M-\sigma_Ks.
\end{equation}
\revred{The minimizers stay a positive distance away from the lateral boundary when $|t-\tau_k^*|\le\delta_K$.}
\end{assumption}

\revred{The two inequalities in \eqref{eq:transversality} express a crossing with a nonzero downward speed. Before the crossing, a time gap $s$ produces a value gap of at least $\sigma_Ks$ above $M$; after the crossing it produces a gap of at least $\sigma_Ks$ below $M$. Consequently a value error of size $e$ can shift the crossing time by at most $e/\sigma_K$. Differentiability of $g_k$ is not required. If $g_k$ is differentiable near $\tau_k^*$ and $g_k'\le-\sigma_K$ there, integration implies \eqref{eq:transversality}. The extra exact time interval is used only to compare crossings that occur slightly earlier or later than $\tau_k^*$; it does not extend the estimates to the quenching time.}

\begin{assumption}[CFL condition and reconstruction positivity]\label{ass:cfl-reconstruction}
\revred{For the same fixed $K$ as in Assumption \ref{ass:finite-stage}, let}
\begin{equation}\label{eq:cfl-h2}
    \tau=\nu h^2,
    \qquad 0<\nu\le\nu_0<\frac12.
\end{equation}
\revred{The reconstruction in \eqref{eq:next-stage-cubic} uses four consecutive old grid nodes, with a one-sided four-node stencil in a boundary cell. For all sufficiently small $h$, its values used to initialize stages $1\le k\le K$ lie in an interval $[a_K,b_K]\subset(0,\infty)$ independent of $h,\tau$.}
\end{assumption}

\revred{The CFL restriction makes the three diffusion weights in the explicit update nonnegative. The positivity requirement keeps the singular reaction $f(s)=-s^{-r}$ in a region where it is Lipschitz: on $[a_K,b_K]$, one may use $|f'(s)|\le ra_K^{-r-1}$. Positive nodal data alone do not guarantee positivity of a cubic interpolant, because some Lagrange weights are negative. The four-node convention also gives a uniform stability bound
$$
    |P_j[Y](x)-P_j[Z](x)|
    \le\Lambda_*\max_{\ell\in S_j}|Y_\ell-Z_\ell|,
$$
where $S_j$ is the chosen stencil and $\Lambda_*$ is independent of $h$ and of the cell. This follows by writing the interpolant in the Lagrange basis and bounding the sum of the absolute values of its four weights.}

\begin{lemma}[One-stage error estimate with non-exact initial data]\label{lem:one-stage-error}
\revblue{Assume the same finite-stage hypotheses as in Assumptions \ref{ass:finite-stage} and \ref{ass:cfl-reconstruction}.}
Let
$$
    \varepsilon_k:=\max_i|U_{i,0}^{(k)}-u^{(k)}(\xi_{k,i},0)|.
$$
\revred{There are constants $C_K>0$ and $\eta_K>0$ such that, if $\varepsilon_k+h^2\le\eta_K$, the Euler sequence can be continued for comparison up to $\tau_k^*+\delta_K$, including the first following grid time, and}
\begin{equation}\label{eq:one-stage-error}
    \max_i|U_{i,n}^{(k)}-u^{(k)}(\xi_{k,i},n\tau)|
    \le C_K(\varepsilon_k+h^2).
\end{equation}
\revred{The same estimate holds for the linear time interpolant on $0\le t\le\tau_k^*+\delta_K$:}
\begin{equation}\label{eq:one-stage-time-interp-error}
    \revblue{\max_i|U_{i,\tau}^{(k)}(t)-u^{(k)}(\xi_{k,i},t)|
    \le C_K(\varepsilon_k+h^2).}
\end{equation}
\revred{The continuation is an auxiliary sequence used in the proof; the algorithm itself stops at its first numerical crossing.}
\end{lemma}
\begin{proof}
{\color{revisionred}
We give the consistency, stability, and positivity arguments separately. Write $u_i^n=u^{(k)}(\xi_{k,i},n\tau)$. Taylor's theorem gives
$$
    \frac{u_i^{n+1}-u_i^n}{\tau}
    =u_t^{(k)}(\xi_{k,i},n\tau)+r^t_{i,n},
    \qquad
    |r^t_{i,n}|\le\frac{\tau}{2}\|u_{tt}^{(k)}\|_\infty.
$$
Expanding the values at $\xi_{k,i}+h$ and $\xi_{k,i}-h$ through order three makes the odd terms cancel. Thus
$$
    \frac{u_{i-1}^n-2u_i^n+u_{i+1}^n}{h^2}
    =u_{\xi\xi}^{(k)}(\xi_{k,i},n\tau)+r^x_{i,n},
    \qquad
    |r^x_{i,n}|\le\frac{h^2}{12}\|u_{\xi\xi\xi\xi}^{(k)}\|_\infty.
$$
The exact PDE therefore implies
$$
    \frac{u_i^{n+1}-u_i^n}{\tau}
    =\delta_x^2u_i^n+f(u_i^n)+R_{i,n},
    \qquad
    |R_{i,n}|\le C_K(\tau+h^2)\le C_Kh^2.
$$
Here $R_{i,n}$ is the defect in the difference equation; its contribution to one update is $\tau R_{i,n}$.

To justify a uniform Lipschitz bound for the numerical sequence as well as the exact one, first use the globally Lipschitz extension
$$
    f_K(s):=-\bigl(\max\{s,m_K/2\}\bigr)^{-r}.
$$
It agrees with $f$ for $s\ge m_K/2$. Its Lipschitz constant is
$$
    L_K=r(m_K/2)^{-r-1}.
$$
The Euler recurrence with $f_K$ is defined throughout the comparison interval. Let $e_{i,n}=U_{i,n}^{(k)}-u_i^n$ and $E_n=\max_i|e_{i,n}|$. Subtracting the exact difference equation from the numerical recurrence gives, with $\nu=\tau/h^2$,
$$
\begin{aligned}
    e_{i,n+1}
    &=(1-2\nu)e_{i,n}+\nu(e_{i-1,n}+e_{i+1,n})\\
    &\quad+\tau\{f_K(U_{i,n}^{(k)})-f_K(u_i^n)\}-\tau R_{i,n}.
\end{aligned}
$$
The boundary errors are zero, because both stages have boundary value $\lambda^{k\alpha}$. Since $1-2\nu\ge0$ by \eqref{eq:cfl-h2}, taking absolute values yields
$$
    E_{n+1}\le(1+L_K\tau)E_n+C_K\tau h^2.
$$
Iteration, with $S_K=\max_{0\le k\le K}(\tau_k^*+2\delta_K)$, gives
$$
    E_n\le e^{L_KS_K}\bigl(\varepsilon_k+C_KS_Kh^2\bigr)
    \qquad\text{for }n\tau\le\tau_k^*+\delta_K+\tau.
$$
Choose $\eta_K$ small enough that this upper bound is at most $m_K/4$, and choose $h$ so that $\tau\le\delta_K$. Then $U_{i,n}^{(k)}\ge3m_K/4>m_K/2$. The cutoff is never used, so this is the original Euler sequence. This proves \eqref{eq:one-stage-error} without assuming in advance that the numerical crossing lies near the exact one.
}

{\color{revisionblue}
Finally, let $0\le t\le\tau_k^*+\delta_K$, choose $t_n\le t\le t_{n+1}$, and set $\theta=(t-t_n)/\tau\in[0,1]$. Both endpoint errors are covered by \eqref{eq:one-stage-error}, because
$$
    t_{n+1}\le\tau_k^*+\delta_K+\tau\le\tau_k^*+2\delta_K.
$$
Add and subtract the linear interpolation of the exact endpoint values to obtain
$$
    U_{i,\tau}^{(k)}(t)-u^{(k)}(\xi_{k,i},t)
    =(1-\theta)e_{i,n}+\theta e_{i,n+1}+r_i(t),
$$
where
$$
    r_i(t):=(1-\theta)u_i^n+\theta u_i^{n+1}-u^{(k)}(\xi_{k,i},t).
$$
The scalar interpolation remainder gives, for some $\zeta_{i,t}\in[t_n,t_{n+1}]$,
$$
    r_i(t)=\frac{\theta(1-\theta)\tau^2}{2}
    u_{tt}^{(k)}(\xi_{k,i},\zeta_{i,t}),
    \qquad
    |r_i(t)|\le\frac{\tau^2}{8}\|u_{tt}^{(k)}\|_\infty.
$$
The factor $1/8$ follows from $\theta(1-\theta)\le1/4$. Taking the maximum over the nodes and using the endpoint estimates yields
$$
\begin{aligned}
    \max_i|U_{i,\tau}^{(k)}(t)-u^{(k)}(\xi_{k,i},t)|
    &\le(1-\theta)E_n+\theta E_{n+1}
    +\frac{\tau^2}{8}\|u_{tt}^{(k)}\|_\infty\\
    &\le C_K(\varepsilon_k+h^2)+C_K\tau^2.
\end{aligned}
$$
Since $\tau=\nu h^2$, the last term is $O(h^4)$ and can be included in $C_K(\varepsilon_k+h^2)$ for small $h$. This proves \eqref{eq:one-stage-time-interp-error}.
}
\end{proof}

\begin{lemma}[Stopping time and next-stage initialization]\label{lem:stage-induction}
\revblue{Assume the same hypotheses as in Lemma \ref{lem:one-stage-error}, including sufficiently small $\varepsilon_k+h^2$.}
\revred{Then}
\begin{equation}\label{eq:stopping-time-error}
    \revred{|\tau_{k,\tau}^*-\tau_k^*|
    \le C_K(\varepsilon_k+h^2).}
\end{equation}
\revred{For $k<K$, the next initial error satisfies}
\begin{equation}\label{eq:next-initial-error}
    \varepsilon_{k+1}\le C_K(\varepsilon_k+h^2).
\end{equation}
\revred{For sufficiently small $\varepsilon_k+h^2$, the next initial nodal minimum is strictly greater than $M$.}
\end{lemma}
\begin{proof}
{\color{revisionred}
\emph{Step 1: replace the continuous minimum by a nodal minimum.}
For $|t-\tau_k^*|\le\delta_K$, let $\xi_*(t)$ be a minimizer and let $\xi_{k,j}$ be a nearest grid node. The interior minimum satisfies $u_\xi^{(k)}(\xi_*(t),t)=0$, and $|\xi_{k,j}-\xi_*(t)|\le h/2$. Taylor's theorem therefore gives
$$
    u^{(k)}(\xi_{k,j},t)-u^{(k)}(\xi_*(t),t)
    =\frac12u_{\xi\xi}^{(k)}(\zeta,t)(\xi_{k,j}-\xi_*(t))^2
    \le\frac18\|u_{\xi\xi}^{(k)}\|_\infty h^2
$$
for some intermediate point $\zeta$. In particular,
\begin{equation}\label{eq:continuous-nodal-minimum}
    0\le\min_i u^{(k)}(\xi_{k,i},t)-g_k(t)\le C_Kh^2.
\end{equation}
For any two finite arrays $a_i,b_i$, one has $|\min_i a_i-\min_i b_i|\le\max_i|a_i-b_i|$. Applying this inequality and Lemma \ref{lem:one-stage-error}, then using \eqref{eq:continuous-nodal-minimum}, yields
\begin{equation}\label{eq:minimum-error-window}
    \sup_{|t-\tau_k^*|\le\delta_K}
    |g_{k,\tau}(t)-g_k(t)|
    \le C_K(\varepsilon_k+h^2).
\end{equation}
The comparison interpolant in this formula is the auxiliary continuation from Lemma \ref{lem:one-stage-error}; its first crossing is the actual stopping time.

\emph{Step 2: locate the first numerical crossing.}
By Assumption \ref{ass:finite-stage},
$0<\delta_K<\min_{0\le k\le K}\tau_k^*$.
For each $k$, the function $g_k$ is continuous, by the joint continuity
of $u^{(k)}$ on a compact space-time cylinder.
The first-crossing property gives $g_k(t)>M$ for $0\le t<\tau_k^*$.
Therefore $g_k-M$ attains a strictly positive minimum on the compact
interval $[0,\tau_k^*-\delta_K]$. Taking the minimum over the finitely
many stages gives
$$
    d_K:=\min_{0\le k\le K}\ \min_{0\le t\le\tau_k^*-\delta_K}
    (g_k(t)-M)>0.
$$
This constant depends only on the exact finite-stage solutions and the
chosen $\delta_K$, and is independent of $h$ and $\tau$.

Fix $t\in[0,\tau_k^*-\delta_K]$. Apply
\eqref{eq:one-stage-time-interp-error} to the auxiliary continuation
from Lemma \ref{lem:one-stage-error}. At every grid node,
$$
    U_{i,\tau}^{(k)}(t)
    \ge u^{(k)}(\xi_{k,i},t)-C_K(\varepsilon_k+h^2)
    \ge g_k(t)-C_K(\varepsilon_k+h^2).
$$
The second inequality holds because $g_k(t)$ is the minimum of the
exact profile over the whole spatial interval, so every exact nodal
value is at least $g_k(t)$. The final lower bound is independent of $i$;
hence it also bounds $g_{k,\tau}(t)=\min_iU_{i,\tau}^{(k)}(t)$ from below.
Decrease the smallness threshold in this lemma so that
$\varepsilon_k+h^2\le d_K/(2C_K)$. Using $g_k(t)\ge M+d_K$, we obtain
$$
\begin{aligned}
    g_{k,\tau}(t)
    &\ge g_k(t)-C_K(\varepsilon_k+h^2)\\
    &\ge M+d_K-\frac{d_K}{2}
     =M+\frac{d_K}{2}>M.
\end{aligned}
$$
Thus no numerical crossing can occur on $[0,\tau_k^*-\delta_K]$.

In the crossing window, choose the error in \eqref{eq:minimum-error-window} less than $\sigma_K\delta_K/2$. Equation \eqref{eq:transversality} then puts $g_{k,\tau}$ above $M$ at the left endpoint and below $M$ at the right endpoint. Continuity gives
$$
    \tau_k^*-\delta_K<\tau_{k,\tau}^*<\tau_k^*+\delta_K,
    \qquad
    g_{k,\tau}(\tau_{k,\tau}^*)=M.
$$
Applying the appropriate side of \eqref{eq:transversality}, according to whether $\tau_{k,\tau}^*$ is earlier or later than $\tau_k^*$, gives
$$
    \sigma_K|\tau_{k,\tau}^*-\tau_k^*|
    \le|g_k(\tau_{k,\tau}^*)-M|
    \le C_K(\varepsilon_k+h^2).
$$
This proves \eqref{eq:stopping-time-error}.

\emph{Step 3: estimate the next initial profile.}
Put $a_i=\lambda^{\beta/2}\xi_{k+1,i}$ and $t_h=\tau_{k,\tau}^*$. Let $j(i)$ select the old cell containing $a_i$, with the one-sided convention of Assumption \ref{ass:cfl-reconstruction} near a boundary. Define $v_\ell=u^{(k)}(\xi_{k,\ell},t_h)$. Linearity of the Lagrange interpolant gives the exact decomposition
$$
\begin{aligned}
    P_{j(i)}[W^{(k)}](a_i)-u^{(k)}(a_i,\tau_k^*)
    &=P_{j(i)}[W^{(k)}-v](a_i)\\
    &\quad+P_{j(i)}[v](a_i)-u^{(k)}(a_i,t_h)\\
    &\quad+u^{(k)}(a_i,t_h)-u^{(k)}(a_i,\tau_k^*).
\end{aligned}
$$
These are, respectively, the propagated nodal error, the spatial interpolation error for a smooth exact profile, and the shift in stopping time. The stability bound for $P_j$, the four-node Lagrange remainder, and the fundamental theorem of calculus imply
\begin{equation}\label{eq:next-error-decomposition}
\begin{aligned}
    |U_{i,0}^{(k+1)}-u^{(k+1)}(\xi_{k+1,i},0)|
    &\le\lambda^\alpha\Lambda_*\|W^{(k)}-v\|_{\ell^\infty}\\
    &\quad+\lambda^\alpha C\|u_{\xi\xi\xi\xi}^{(k)}\|_\infty h^4\\
    &\quad+\lambda^\alpha\|u_t^{(k)}\|_\infty|t_h-\tau_k^*|.
\end{aligned}
\end{equation}
For clarity, the Lagrange remainder is
$$
    P_{j(i)}[v](a_i)-u^{(k)}(a_i,t_h)
    =-\frac{u_{\xi\xi\xi\xi}^{(k)}(\zeta_i,t_h)}{4!}
    \prod_{\ell\in S_{j(i)}}(a_i-\xi_{k,\ell});
$$
each factor is at most $3h$ in absolute value. At a boundary node the reconstruction and exact boundary data coincide, so the corresponding error is zero. Also, $W_\ell^{(k)}=U_{\ell,\tau}^{(k)}(t_h)$, and hence
\begin{equation}\label{eq:W-stage-error}
    \|W^{(k)}-v\|_{\ell^\infty}\le C_K(\varepsilon_k+h^2)
\end{equation}
by \eqref{eq:one-stage-time-interp-error}. Substituting \eqref{eq:W-stage-error} and \eqref{eq:stopping-time-error} into \eqref{eq:next-error-decomposition} proves \eqref{eq:next-initial-error}. The exact next initial minimum is $\lambda^\alpha M>M$. By \eqref{eq:next-initial-error},
$\varepsilon_{k+1}\le C_K(\varepsilon_k+h^2)$.
Thus we can decrease the threshold for $\varepsilon_k+h^2$ so that
$$
\varepsilon_k+h^2<\frac{(\lambda^\alpha-1)M}{2C_K},
\qquad
\varepsilon_{k+1}<\frac{(\lambda^\alpha-1)M}{2}.
$$
For exact initial sampling, $\varepsilon_0=0$, finite iteration of
\eqref{eq:next-initial-error} gives $\varepsilon_k\le C_Kh^2$, as detailed
in the proof of Theorem \ref{thm:convergence} below. Therefore
$$
    \min_iU_{i,0}^{(k+1)}
    \ge\lambda^\alpha M-\varepsilon_{k+1}>M.
$$
Thus the next stage starts strictly above the stopping level.

Here and throughout this subsection, $C_K$ is a generic constant: it may be increased from line to line. It includes the finitely many derivative bounds, $1/\sigma_K$, the fixed factor $\lambda^\alpha$, and the interpolation stability constant. None depends on $h$ or $\tau$. No bound uniform as $K\to\infty$ is asserted.
}
\end{proof}

\begin{theorem}[Finite-stage convergence of the rescaling algorithm]\label{thm:convergence}
\revblue{Fix $K\in\mathbb N$. Assume Assumptions \ref{ass:finite-stage} and \ref{ass:cfl-reconstruction}, the same hypotheses used in Lemmas \ref{lem:one-stage-error} and \ref{lem:stage-induction}, and sample the initial data exactly: $U_{i,0}^{(0)}=u_0(x_i)$.}
\revred{There exist $h_K>0$ and $C_K>0$ such that, for $0<h\le h_K$ and $\tau=\nu h^2$,}
\begin{equation}\label{eq:finite-stage-convergence}
    \revred{\max_{0\le k\le K}\ \max_{0\le n\le n_k}\ \max_i
    |U_{i,n}^{(k)}-u^{(k)}(\xi_{k,i},n\tau)|
    \le C_Kh^2.}
\end{equation}
\revred{The grid time $n_k\tau$ is the first full Euler time at which the nodal minimum is at or below $M$ and is included in this estimate, although the interpolated stopping time can be earlier in that step. These endpoints lie inside the enlarged exact intervals when $h\le h_K$.}
Moreover,
\begin{equation}\label{eq:muK-convergence}
    \revred{\max_{0\le k\le K}|\tau_{k,\tau}^*-\tau_k^*|
    +\max_{0\le k\le K}|\mu_{k,h,\tau}-\mu_k|
    \le C_Kh^2.}
\end{equation}
\revred{At the physical stopping times the numerical nodal stopping level is
$$
    \min_i\lambda^{-\alpha k}W_i^{(k)}=\lambda^{-\alpha k}M,
$$
and it satisfies}
\begin{equation}\label{eq:physical-stopping-error}
    \revred{\left|
    \lambda^{-\alpha k}M-\min_{x\in[-1,1]}u(x,\mu_{k,h,\tau})
    \right|\le C_Kh^2,\qquad 0\le k\le K.}
\end{equation}
\revred{The stopping level is the minimum over nodes; equality with the minimum of a continuous cubic interpolant is not needed.}
\end{theorem}
\begin{proof}
{\color{revisionred}
Exact sampling gives $\varepsilon_0=0$. Lemma \ref{lem:stage-induction} gives $\varepsilon_{k+1}\le C_K(\varepsilon_k+h^2)$, and iteration over the finite set of stages gives $\varepsilon_k\le C_Kh^2$ for $0\le k\le K$. More explicitly, if $C_*$ is a common constant in that lemma, one may take $B_0=0$ and $B_{k+1}=C_*(B_k+1)$, so that $\varepsilon_k\le B_kh^2$. Choosing $h_K$ so that $(1+\max_{k\le K}B_k)h_K^2$ is below the smallness thresholds of Lemmas \ref{lem:one-stage-error} and \ref{lem:stage-induction} makes every induction step valid.

Lemma \ref{lem:one-stage-error} now proves \eqref{eq:finite-stage-convergence}. By Lemma \ref{lem:stage-induction},
$$
    n_k\tau\le\tau_{k,\tau}^*+\tau
    \le\tau_k^*+C_Kh^2+\tau<\tau_k^*+\delta_K
$$
for small $h$. This justifies including the endpoint of the crossing step. The estimate for the stopping times, together with
$$
    |\mu_{k,h,\tau}-\mu_k|
    \le\sum_{j=0}^k\lambda^{\beta j}|\tau_{j,\tau}^*-\tau_j^*|
    \le C_Kh^2\sum_{j=0}^k\lambda^{\beta j}
$$
prove \eqref{eq:muK-convergence}.
}

{\color{revisionblue}
For \eqref{eq:physical-stopping-error}, write $m(t)=\min_{x\in[-1,1]}u(x,t)$. The exact finite-stage intervals cover a compact physical time interval slightly beyond $\mu_K$, on which $u_t$ is bounded. Indeed, the exact scaling relation gives
$$
    u^{(j)}(\xi,\eta)=\lambda^{j\alpha}
    u\bigl(\lambda^{\beta j/2}\xi,\mu_{j-1}+\lambda^{\beta j}\eta\bigr),
$$
so one can bound the physical time derivative by
$$
    L_K^{\mathrm{phys}}
    :=\max_{0\le j\le K}\lambda^{-j(\alpha+\beta)}
    \|u_t^{(j)}\|_{L^\infty(Q_{j,K})}<\infty.
$$
This constant depends on the finitely many bounds in Assumption \ref{ass:finite-stage} and the fixed scaling parameters. The numerical times also belong to this interval for sufficiently small $h$, by \eqref{eq:muK-convergence}. For any two times $s,t$ in it, the fundamental theorem of calculus gives
$$
    u(x,t)-u(x,s)=\int_s^t u_t(x,\zeta)\,d\zeta.
$$
Taking absolute values and the supremum over $x$ yields
$$
    \|u(\cdot,t)-u(\cdot,s)\|_\infty
    \le\sup_{x,\zeta}|u_t(x,\zeta)|\,|t-s|
    \le L_K^{\mathrm{phys}}|t-s|.
$$
To pass to the minima, observe that every $x\in[-1,1]$ satisfies
$$
\begin{aligned}
    u(x,t)&\ge u(x,s)-\|u(\cdot,t)-u(\cdot,s)\|_\infty\\
    &\ge m(s)-\|u(\cdot,t)-u(\cdot,s)\|_\infty.
\end{aligned}
$$
Taking the minimum over $x$ and then interchanging $s$ and $t$ proves
$$
    |m(t)-m(s)|\le\|u(\cdot,t)-u(\cdot,s)\|_\infty
    \le L_K^{\mathrm{phys}}|t-s|.
$$
}
\revred{At the exact accumulated stopping time, $m(\mu_k)=\lambda^{-\alpha k}M$. Consequently
$$
    |m(\mu_{k,h,\tau})-\lambda^{-\alpha k}M|
    \le L_K^{\mathrm{phys}}|\mu_{k,h,\tau}-\mu_k|
    \le C_Kh^2,
$$
which proves \eqref{eq:physical-stopping-error}.}
\end{proof}

\subsection{Convergence of the numerical quenching time}
\revred{The finite-stage hypotheses are exactly Assumptions \ref{ass:finite-stage} and \ref{ass:cfl-reconstruction}, as in Theorem \ref{thm:convergence}. In this subsection they are required for every fixed finite number of stages. Passing to infinitely many stages needs an additional uniform tail hypothesis, which we state separately.}

\begin{assumption}[\revred{Admissible tail stages and comparison}]\label{ass:tail-comparison}
{\color{revisionred}
There exist an integer $K_{\rm tail}\ge1$, a mesh threshold $h_{\rm tail}>0$, and a constant $c_0\in(0,1]$ with the following properties. Use the parameter conditions \eqref{parameters}, the step restriction in Assumption \ref{ass:cfl-reconstruction}, and $h_{\rm tail}$ small enough that $\tau<M^{r+1}$. For $0<h\le h_{\rm tail}$ the numerical algorithm has admissible stages for all $k\ge K_{\rm tail}$: each reconstructed initial profile satisfies $\min_iU_{i,0}^{(k)}>M$, and the Euler values used through the crossing step are positive. At every pre-crossing time level $0\le n<n_k$, there is an interior nodal minimizer $i_{k,n}$ satisfying
\begin{equation}\label{eq:tail-comparison-ineq}
    \delta_x^2U_{i_{k,n},n}^{(k)}
    \le(1-c_0)(U_{i_{k,n},n}^{(k)})^{-r}.
\end{equation}
The constants $c_0$ and $h_{\rm tail}$ are independent of $k$.
}
\end{assumption}

\revred{Here \eqref{eq:tail-comparison-ineq} is an additional comparison assumption; it is not a consequence of the finite-stage error estimate. It says that diffusion cancels at most a fraction $1-c_0$ of the singular reaction at a nodal minimum. If the hypotheses of Proposition \ref{prop:accumulated-defect}, including its defect-absorption condition, hold uniformly on a region containing those minima, its conclusion \eqref{ineq:absorbed-cprime} gives this inequality with $c_0=c/2$. This is only one conditional sufficient route. In particular, the constant $c_0$ is the comparison margin defined here, and is not a new discretization parameter. The computations in Section \ref{sec:numerics} check the margin on finitely many stages; they do not establish the infinite-stage assumption.}

\begin{prp}[Tail bound]\label{prop:pde-tail-bound}
\revblue{Use the same aligned-grid algorithm and step restrictions as in Theorem \ref{thm:convergence}, and assume Assumption \ref{ass:tail-comparison}.}
\revred{Then, for every $K_0\ge K_{\rm tail}$ and all sufficiently small $h$, there is a constant $C>0$, independent of $k,K_0,h,\tau$, such that}
\begin{equation}\label{eq:pde-tail-bound}
    T_{h,\tau}-\mu_{K_0-1,h,\tau}
    \le\lambda^{\beta K_0}
    \left(
        \frac{(\lambda^\alpha M)^{r+1}}{(r+1)c_0}
        +C\frac{\tau}{1-\lambda^\beta}
    \right).
\end{equation}
\revred{The regularity and transversality hypotheses of Theorem \ref{thm:convergence} are not needed for this tail estimate itself.}
\end{prp}
\begin{proof}
{\color{revisionred}
Put $m_n^{(k)}=\min_iU_{i,n}^{(k)}$. At the minimizer supplied by Assumption \ref{ass:tail-comparison}, the Euler scheme gives
$$
\begin{aligned}
    m_{n+1}^{(k)}
    &\le U_{i_{k,n},n+1}^{(k)}\\
    &=m_n^{(k)}+\tau\left\{\delta_x^2U_{i_{k,n},n}^{(k)}
    -(m_n^{(k)})^{-r}\right\}\\
    &\le m_n^{(k)}-c_0\tau(m_n^{(k)})^{-r}.
\end{aligned}
$$
This computation explains why $c_0$ occurs in the denominator of the tail bound.

We next bound the initial minimum of every stage $k\ge1$. Let $j_{k-1}$ satisfy $W_{j_{k-1}}^{(k-1)}=M$, as guaranteed by \eqref{eq:min-W-equals-M}. Since $q=\lambda^{-\beta/2}$ is an integer, the next grid contains $i=qj_{k-1}$ and
$$
    \lambda^{\beta/2}\xi_{k,i}=\xi_{k-1,j_{k-1}}.
$$
At this retained old node $i_*=qj_{k-1}$, the interpolation property gives $U_{i_*,0}^{(k)}=\lambda^\alpha M$. More explicitly,
$$
\min_{-I_{k-1}\le j\le I_{k-1}}U_{qj,0}^{(k)}=\lambda^\alpha M,
\qquad
\min_{-I_k\le i\le I_k}U_{i,0}^{(k)}\le U_{i_*,0}^{(k)}=\lambda^\alpha M.
$$
The first minimum is taken only over retained old nodes. Cubic interpolation may produce smaller values at newly inserted nodes, so equality for the full-grid minimum is not required. Together with the strict restart condition, this gives
\begin{equation}\label{eq:tail-initial-min}
    M<m_0^{(k)}\le a_\lambda,
    \qquad a_\lambda:=\lambda^\alpha M.
\end{equation}
The restriction $K_{\rm tail}\ge1$ ensures that the first, unrescaled initial datum is not used in this estimate.
}

{\color{revisionblue}
For clarity, the exact solution of $w'=-c_0w^{-r}$ with $w(0)=a_\lambda$ reaches $M$ at the scalar reference time
$$
    t_M(a_\lambda;c_0)
    =\int_M^{a_\lambda}\frac{s^r}{c_0}\,ds
    =\frac{a_\lambda^{r+1}-M^{r+1}}{(r+1)c_0}.
$$
This identity follows by separating variables and explains the main term in the local crossing estimate.

Consider the scalar Euler sequence
}
{\color{revisionred}
$$
    z_{n+1}=z_n-c_0\tau z_n^{-r},
    \qquad z_0=a_\lambda.
$$
The map $s\mapsto s-c_0\tau s^{-r}$ is increasing on $(0,\infty)$, because its derivative is $1+c_0r\tau s^{-r-1}>0$. Induction therefore gives $m_n^{(k)}\le z_n$ as long as both sequences are above $M$. Consequently the PDE minimum crosses $M$ no later than the first full Euler time level at which $z_n\le M$.
}
{\color{revisionblue}
Let $\widehat\tau_\tau$ denote the scalar linearly interpolated crossing time. Lemma \ref{lem:ode-local-crossing}, with $b=c_0$ and $A=a_\lambda$, gives
$$
    |\widehat\tau_\tau-t_M(a_\lambda;c_0)|\le C\tau.
$$
If $N$ is the first index for which $z_N\le M$, then $N\tau\le\widehat\tau_\tau+\tau$. Combining the comparison with this possible final full step, and incorporating its contribution into $C$, yields
\begin{equation}\label{eq:tail-local-bound}
    \tau_{k,\tau}^*
    \le\frac{a_\lambda^{r+1}-M^{r+1}}{(r+1)c_0}+C\tau.
\end{equation}
Assumption \ref{ass:tail-comparison} supplies a common $c_0$ for all tail stages and meshes. The parameters $a_\lambda,M,r,c_0$ in the scalar crossing estimate are fixed, so its step threshold and the constant $C$ can be chosen uniformly in $k$, $h$, and $\tau$.
}

{\color{revisionred}
Multiplying \eqref{eq:tail-local-bound} by $\lambda^{\beta k}$ and summing first over finitely many stages and then passing to the increasing limit yields
$$
\begin{aligned}
    T_{h,\tau}-\mu_{K_0-1,h,\tau}
    &\le\frac{\lambda^{\beta K_0}}{1-\lambda^\beta}
    \left\{\frac{a_\lambda^{r+1}-M^{r+1}}{(r+1)c_0}
    +C\tau\right\}\\
    &=\lambda^{\beta K_0}
    \left\{\frac{a_\lambda^{r+1}}{(r+1)c_0}
    +\frac{C\tau}{1-\lambda^\beta}\right\}.
\end{aligned}
$$
The last equality follows from
$$
    a_\lambda^{r+1}=\lambda^{-\beta}M^{r+1},
    \qquad
    \frac{a_\lambda^{r+1}-M^{r+1}}{1-\lambda^\beta}
    =a_\lambda^{r+1}.
$$
This proves both finiteness of the numerical tail and \eqref{eq:pde-tail-bound}.
}
\end{proof}

\begin{theorem}[Conditional convergence of the numerical quenching time]\label{thm:pde-quenching-time-convergence}
\revblue{Assume exactly the finite-stage hypotheses of Theorem \ref{thm:convergence}, namely Assumptions \ref{ass:finite-stage} and \ref{ass:cfl-reconstruction}, for every fixed $K$. In addition, assume the admissibility and uniform tail comparison of Assumption \ref{ass:tail-comparison}.}
Then
\begin{equation}\label{eq:pde-time-convergence}
    \lim_{h\to0,\ \tau=\nu h^2}T_{h,\tau}=T.
\end{equation}
\end{theorem}
\begin{proof}
{\color{revisionred}
At the exact stopping time $\mu_k$, the physical minimum reaches the stopping level $\lambda^{-\alpha k}M$. The numerical solution reaches the same nodal stopping level at the numerical stopping time $\mu_{k,h,\tau}$. Theorem \ref{thm:convergence} gives convergence of these stopping times for each fixed $k$.

First, $\mu_k\uparrow T$. Indeed, the exact first crossings occur before quenching and their levels $\lambda^{-\alpha k}M$ decrease to zero. If their increasing limit were $\bar T<T$, continuity and positivity on $[-1,1]\times[0,\bar T]$ would imply a positive lower bound for the minimum, contradicting $\min_xu(x,\mu_k)=\lambda^{-\alpha k}M\to0$.

For each fixed $K$, $T_{h,\tau}\ge\mu_{K,h,\tau}$. Hence the stopping-time estimate \eqref{eq:muK-convergence} gives
$$
    \liminf_{h\to0}T_{h,\tau}\ge\mu_K.
$$
Letting $K\to\infty$ proves $\liminf_{h\to0}T_{h,\tau}\ge T$.

For the upper bound, fix $K_0\ge K_{\rm tail}$. Proposition \ref{prop:pde-tail-bound} and the same finite-stage stopping-time estimate give
$$
    \limsup_{h\to0}T_{h,\tau}
    \le\mu_{K_0-1}+\lambda^{\beta K_0}
    \frac{(\lambda^\alpha M)^{r+1}}{(r+1)c_0}.
$$
Here $\tau=\nu h^2\to0$, while $K_0$ is held fixed. Now let $K_0\to\infty$. Since $0<\lambda^\beta<1$ and $\mu_{K_0-1}\uparrow T$, the right-hand side tends to $T$. This proves the upper bound and completes the proof.
}
\end{proof}

\section{Numerical implementation and results}\label{sec:numerics}
This section reports aligned-grid numerical experiments for the rescaling algorithm applied to the one-dimensional problem
\begin{equation}\label{eq:numerical-model}
    u_t=u_{xx}-u^{-r},\qquad -1<x<1,
\end{equation}
with boundary condition $u(\pm1,t)=1$.  The finite-difference evolution is the explicit Euler scheme \eqref{scheme}.  At the end of each rescaled stage, the stopping time is located by the node-wise interpolation rule \eqref{eq:nodal-interpolated-minimum} to \eqref{eq:W-stopping-profile}, and the next initial profile is reconstructed in space by the local cubic interpolation \eqref{eq:next-stage-cubic}.  Thus the numerical implementation follows the same two interpolation rules used in the preceding analysis: linear interpolation only in time, and cubic interpolation in space at the rescaling step.

The threshold used by the rescaling algorithm is denoted by $M$, while the minimum of the initial data is denoted by $M_0$ in order to avoid the ambiguity $\min u_0=M$ at $t=0$.  We take
\begin{equation}\label{eq:numerical-parameters}
    r=6,
    \qquad
    M=0.5,
    \qquad
    M_0=0.7,
    \qquad
    \beta=2,
    \qquad
    \lambda=\frac12,
    \qquad
    \alpha=-\frac27.
\end{equation}
Then $q=\lambda^{-\beta/2}=2\in\mathbb N$, so the computation matches the aligned-grid assumption \eqref{eq:aligned-q}.  The initial profile is
\begin{equation}\label{eq:numerical-initial-data}
    u_0(x)=(1-M_0)x^2+M_0.
\end{equation}
Thus $\min_xu_0(x)=M_0>M$, and the first rescaling is not performed at $t=0$.  For each base mesh $I$, the initial grid is $x_i=i/I$ on $[-1,1]$, so $h=1/I$.  The time step is chosen as
\begin{equation}\label{eq:numerical-cfl}
    \tau=0.25h^2,
\end{equation}
{\color{revisionblue}
which satisfies the explicit CFL condition. The computation is run up to $K_{\max}=8$. For the remaining tail, we use the crossing time of the scalar equation $w'=-w^{-r}$ from $w(0)=\lambda^\alpha M$ to $M$ as a reference. Separation of variables and the scaling relation $\alpha(r+1)=-\beta$ give
$$
\int_M^{\lambda^\alpha M}s^r\,ds
=\frac{(\lambda^\alpha M)^{r+1}-M^{r+1}}{r+1}
=\frac{M^{r+1}}{r+1}(\lambda^{-\beta}-1).
$$
We denote this scalar reference time by
\begin{equation}\label{eq:tau-infty-numerical}
\tau_\infty=\frac{M^{r+1}}{r+1}(\lambda^{-\beta}-1).
\end{equation}
The value $\tau_\infty$ is the exact crossing time for this scalar ODE. Its use for PDE stages corresponds to a regime in which diffusion near the relevant minimum is small compared with the singular reaction. Establishing convergence of the PDE stage times to this value requires additional asymptotic control of the rescaled profiles and their diffusion terms. For scalar Euler stages, Lemma~\ref{lem:ode-local-crossing} gives an $O(\tau)$ difference from this reference time. Convergence to the exact reference therefore also requires control of temporal discretization in the chosen limiting regime. For the reported numerical estimator, we approximate each uncomputed local stage time by $\tau_\infty$. Summing the associated physical-time weights gives
$$
\sum_{k=K_{\max}+1}^\infty\lambda^{\beta k}\tau_\infty
=\frac{\lambda^{\beta(K_{\max}+1)}}{1-\lambda^\beta}\tau_\infty.
$$
The reported quantity is the tail-corrected estimator of the quenching time
\begin{equation}\label{eq:numerical-tail-corrected-time}
T_{h,\tau}^{(K_{\max})}
=\sum_{k=0}^{K_{\max}}\lambda^{\beta k}\tau_{k,\tau}^*
+\frac{\lambda^{\beta(K_{\max}+1)}}{1-\lambda^\beta}\tau_\infty.
\end{equation}
}
For three successive meshes $I,2I,4I$, the observed order is computed by
\begin{equation}\label{eq:observed-order}
    p_{4I}
    =\log_2\frac{|T_I^{(K_{\max})}-T_{2I}^{(K_{\max})}|}
                   {|T_{2I}^{(K_{\max})}-T_{4I}^{(K_{\max})}|},
    \qquad I=20,40.
\end{equation}
Thus the order displayed in the row $4I$ uses the three resolutions $I,2I,4I$.

\begin{table}[htbp]
\centering
\begin{tabular}{ccccc}
\toprule
$I$ & $h$ & $\tau$ & $T_{h,\tau}^{(K_{\max})}$ & observed order \\
\midrule
20 & 0.05000 & $6.250\times10^{-4}$ & 0.0133016787 & -- \\
40 & 0.02500 & $1.563\times10^{-4}$ & 0.0128155098 & -- \\
80 & 0.01250 & $3.906\times10^{-5}$ & 0.0126846130 & 1.893 \\
160 & 0.00625 & $9.766\times10^{-6}$ & 0.0126510836 & 1.965 \\
\bottomrule
\end{tabular}
\caption{Tail-corrected estimates of the quenching time for the aligned-grid test.  The order is estimated from the successive-difference formula \eqref{eq:observed-order}.}
\label{tab:computed-times}
\end{table}

{\color{revisionred}
Table~\ref{tab:computed-times} shows that the tail-corrected estimates stabilize under mesh refinement and that the observed order is close to two for the finest refinements. This is an empirical observation for the stated data and the choice $\tau=O(h^2)$.

The hypotheses of Theorem~\ref{thm:convergence} need separate verification before that theorem can be applied to a particular example. In fact, the initial profile \eqref{eq:numerical-initial-data} does not satisfy the boundary compatibility required for $C^{4,2}$ regularity on the closed space-time cylinder starting at $t=0$: the constant boundary data require $u_t(\pm1,0)=0$, whereas
$$
u_0''(\pm1)-u_0(\pm1)^{-r}=2(1-M_0)-1=-0.4.
$$
Thus the present finite-stage theorem, with its stated up-to-boundary regularity assumption, does not directly cover this numerical example. A convergence proof for these data would need an analysis of the initial boundary layer under weaker regularity assumptions. The computations also do not prove the uniform infinite-stage tail hypothesis.
}

To check the tail comparison condition a posteriori, we also compute at the nodal minimizer $i_{k,n}$ the ratio
\begin{equation}\label{eq:numerical-tail-ratio}
    R_{k,n}
    =\frac{\delta_x^2U_{i_{k,n},n}^{(k)}}{(U_{i_{k,n},n}^{(k)})^{-r}}.
\end{equation}
Assumption \ref{ass:tail-comparison} requires $R_{k,n}\le1-c_0$ in the tail.  Table \ref{tab:tail-ratio} reports the maximum of this ratio for $4\le k\le8$ over all time levels before stopping.  These values do not prove the tail hypothesis, but they show that the numerical run is consistent with a uniform choice such as $c_0=0.95$: in the table, $\max R_{k,n}<0.05$, hence $R_{k,n}\le1-c_0$ for $c_0=0.95$ on the reported tail stages.

\begin{table}[htbp]
\centering
\begin{tabular}{ccc}
\toprule
$I$ & $\max_{4\le k\le8,n} R_{k,n}$ & corresponding lower margin $1-\max R_{k,n}$ \\
\midrule
40  & 0.04324 & 0.95676 \\
80  & 0.04610 & 0.95390 \\
160 & 0.04687 & 0.95313 \\
\bottomrule
\end{tabular}
\caption{A posteriori check of the tail comparison ratio at nodal minima.}
\label{tab:tail-ratio}
\end{table}

\begin{figure}[htbp]
  \centering
  \includegraphics[width=0.78\textwidth]{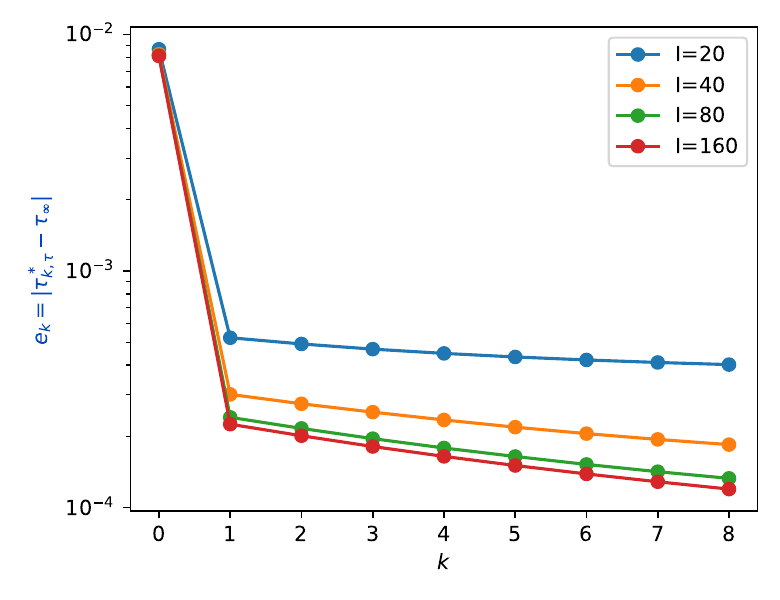}
  \caption{\revblue{Absolute differences $e_k=|\tau_{k,\tau}^*-\tau_\infty|$ from the scalar reference time for $I=20,40,80,160$ and $k=0,\ldots,8$. The dependence on $h$ is suppressed in the stopping-time notation.}}
  \label{fig:error}
\end{figure}

\begin{figure}[htbp]
  \centering
  \includegraphics[width=0.78\textwidth]{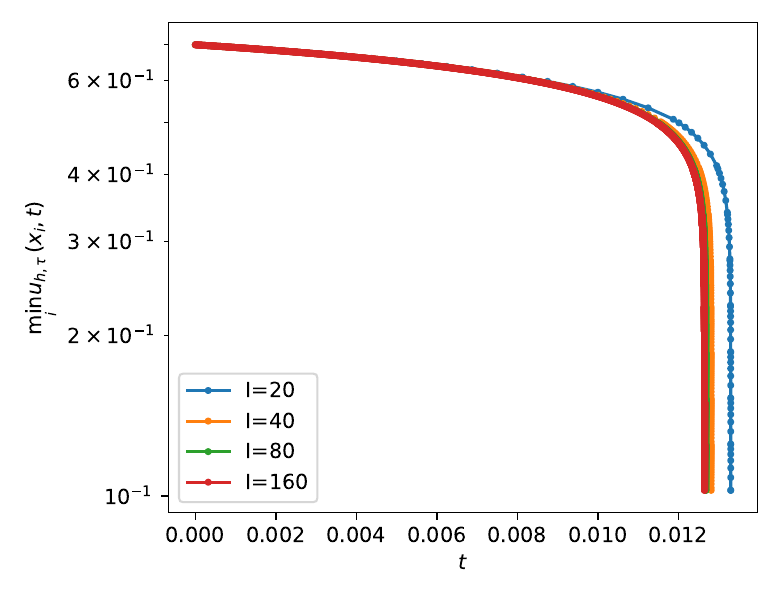}
  \caption{Evolution of the physical nodal minimum for $I=20,40,80,160$.  Full Euler values before each crossing and the interpolated stopping points are shown.}
  \label{fig:min-u}
\end{figure}

\begin{figure}[htbp]
  \centering
  \includegraphics[width=0.78\textwidth]{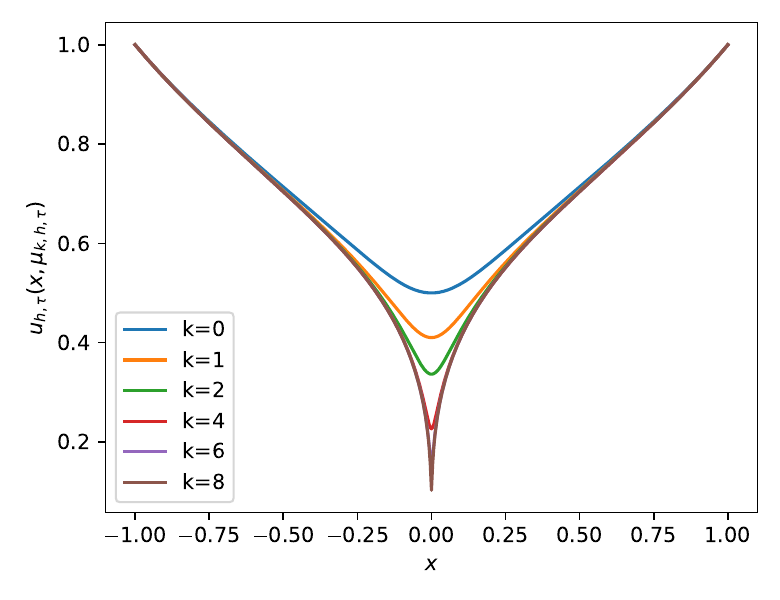}
  \caption{Representative profiles from the finest computation at the physical stopping instants $\mu_{k,h,\tau}$, with $k=0,1,2,4,6,8$.}
  \label{fig:solution}
\end{figure}

Figure~\ref{fig:error} plots
{\color{revisionred}
\begin{equation}\label{eq:error-tau-star}
e_k=|\tau_{k,\tau}^*-\tau_\infty|,
\end{equation}
}
\revblue{Here $\tau_\infty$ is the scalar reference time in \eqref{eq:tau-infty-numerical}. For the reported meshes and stages, these deviations decrease as $k$ increases and under spatial refinement.}
Figure~\ref{fig:min-u} shows the decay of the physical nodal minimum toward zero for the same meshes as in Table~\ref{tab:computed-times}. Figure~\ref{fig:solution} shows representative physical profiles from the finest computation.

\clearpage
\section*{Acknowledgments}
Maissâ Boughrara Boughrara gratefully acknowledges the Fondation Sciences Mathématiques de Paris for funding a two-month research visit at Musashino University, which contributed to the development of this work. The authors thank Professor Tetsuji Tokihiro for valuable ideas and helpful discussions related to this project. The research of Maissâ Boughrara was funded in whole or in part by the Austrian Science Fund (FWF) 10.55776/PAT5825523.

\appendix
\section{\texorpdfstring{Proof of estimate \eqref{MaxPrin J}}{Proof of the discrete defect estimate}}\label{proof max prin J}

In this section, we give the detailed computations for estimate \eqref{MaxPrin J}.
\begin{proof}[Proof of \eqref{MaxPrin J}]
We work at an interior node where all required difference stencils are
defined and the Euler equation holds at neighboring nodes. The values involved
are positive and $0\le c\le1$.
We denote by $V^{(k)}_{i,n}=\delta_t U^{(k)}_{i,n}$.
With \eqref{def J}, Taylor expansion and Proposition
\ref{prop:pde-positivity}, we have that
\begin{equation}\label{delta t J upper bound}
\begin{aligned}
\delta_tJ_{i,n}^{(k)}
&=\delta_t\delta_x^2 U_{i,n}^{(k)}
 -(1-c)\delta_t((U_{i,n}^{(k)})^{-r})\\
&=\delta_x^2V_{i,n}^{(k)}
 -(1-c)\delta_t((U_{i,n}^{(k)})^{-r})\\
&=\delta_x^2V_{i,n}^{(k)}
 +r(1-c)(U_{i,n}^{(k)})^{-r-1}V_{i,n}^{(k)}\\
&\quad -(1-c)\frac{r(r+1)}{2}\xi_{i,n}^{-r-2}
 \frac{(U_{i,n+1}^{(k)}-U_{i,n}^{(k)})^2}{\tau}\\
&\leq \delta_x^2V_{i,n}^{(k)}
 +r(1-c)(U_{i,n}^{(k)})^{-r-1}V_{i,n}^{(k)}\\
&=\delta_x^2V_{i,n}^{(k)}
 +r(1-c)(U_{i,n}^{(k)})^{-r-1}J_{i,n}^{(k)}\\
&\quad -rc(1-c)(U_{i,n}^{(k)})^{-2r-1},
\end{aligned}
\end{equation}
where $\xi_{i,n}$ lies between $U_{i,n+1}^{(k)}$ and $U_{i,n}^{(k)}$.
On the other hand, we have that
$$
\delta_x^2J_{i,n}^{(k)}
=\delta_x^2V_{i,n}^{(k)}+c\delta_x^2((U_{i,n}^{(k)})^{-r}).
$$
Again with Taylor expansion and Proposition \ref{prop:pde-positivity}, we have that
$$
\begin{aligned}
\delta_x^2((U_{i,n}^{(k)})^{-r})
&=-r(U_{i,n}^{(k)})^{-r-1}\delta_x^2U_{i,n}^{(k)}\\
&\quad+\frac{r(r+1)}{2h^2}
 (\xi_{i,n}^{(k)})^{-r-2}
 (U_{i+1,n}^{(k)}-U_{i,n}^{(k)})^2\\
&\quad+\frac{r(r+1)}{2h^2}
 (\bar\xi_{i,n}^{(k)})^{-r-2}
 (U_{i-1,n}^{(k)}-U_{i,n}^{(k)})^2\\
&\geq-r(U_{i,n}^{(k)})^{-r-1}\delta_x^2U_{i,n}^{(k)},
\end{aligned}
$$
for some $\xi^{(k)}_{i,n},\bar\xi^{(k)}_{i,n}$ between
$U^{(k)}_{i+1,n},U^{(k)}_{i,n}$ and
$U^{(k)}_{i-1,n},U^{(k)}_{i,n}$, respectively.
Therefore, with \eqref{scheme} and \eqref{def J}, we have
$$
\begin{aligned}
\delta_x^2J_{i,n}^{(k)}
&\geq\delta_x^2V_{i,n}^{(k)}
 -cr(U_{i,n}^{(k)})^{-r-1}\delta_x^2U_{i,n}^{(k)}\\
&=\delta_x^2V_{i,n}^{(k)}
 -cr(U_{i,n}^{(k)})^{-r-1}J_{i,n}^{(k)}\\
&\quad-cr(1-c)(U_{i,n}^{(k)})^{-2r-1}.
\end{aligned}
$$
Then, with \eqref{delta t J upper bound}, we obtain
$$
\delta_tJ_{i,n}^{(k)}-\delta_x^2J_{i,n}^{(k)}
\leq r(U_{i,n}^{(k)})^{-r-1}J_{i,n}^{(k)}.
$$
This concludes the proof of estimate \eqref{MaxPrin J}.
\end{proof}

\bibliographystyle{alpha}
\bibliography{main.bbl}
\end{document}